\documentclass[11pt,a4paper]{article}
\usepackage[T1]{fontenc}
\usepackage[utf8]{inputenc}
\usepackage{lmodern}
\usepackage[english]{babel}
\usepackage{amsmath,amssymb,amsthm,mathtools}
\usepackage[margin=26mm]{geometry}
\usepackage{booktabs,array,enumitem,graphicx,microtype}
\usepackage[font=small,labelfont=bf,skip=8pt]{caption}
\usepackage[hidelinks]{hyperref}
\hypersetup{pdftitle={Novel upper bounds for Borsuk numbers from contact points and spherical partitions},
  pdfkeywords={Borsuk problem, spherical partitions, contact points, computer-assisted proof}}

\newtheorem{mainthm}{Theorem}
\newtheorem{theorem}{Theorem}[section]
\newtheorem{lemma}[theorem]{Lemma}
\newtheorem{proposition}[theorem]{Proposition}
\newtheorem{corollary}[theorem]{Corollary}
\theoremstyle{definition}
\newtheorem{definition}[theorem]{Definition}
\theoremstyle{remark}
\newtheorem{remark}[theorem]{Remark}

\numberwithin{equation}{section}
\renewcommand{\arraystretch}{1.1}

\newcommand{\R}{\mathbb R}
\newcommand{\Q}{\mathbb Q}
\newcommand{\Z}{\mathbb Z}
\newcommand{\ip}[2]{\langle #1,#2\rangle}
\newcommand{\diam}{\operatorname{diam}}
\newcommand{\conv}{\operatorname{conv}}
\newcommand{\cone}{\operatorname{cone}}
\newcommand{\aff}{\operatorname{aff}}
\newcommand{\vertset}{\operatorname{vert}}
\newcommand{\eps}{\varepsilon}
\newcommand{\muhat}{\widehat\mu}
\newcommand{\uround}{\mathrm{u}}
\newcommand{\repo}{\url{https://github.com/Alexandr-Tolmachev/borsuk-bounds}}

\title{Novel upper bounds for Borsuk numbers from \\ contact points and spherical partitions}
\author{Alexander Tolmachev\footnote{Moscow Institute of Physics and Technology; email: tolmachev.ad@phystech.edu}}

\date{}

\begin{document}
\maketitle

\begin{abstract}
Let $b(n)$ be the least $m$ such that every bounded set in $\R^n$ of
positive diameter can be partitioned into $m$ sets of smaller diameter. We
prove that
\[
 b(4)\le6,\qquad b(5)\le12,\qquad b(6)\le31,\qquad b(7)\le45,\qquad b(8)\le89,
\]
and that
\[
 b(n)\le1+2^{n-1}\left(\tfrac{39}{64}\right)^{\lfloor(n-1)/8\rfloor}
 \qquad\text{for every } n\ge9 ;
\]
for instance, $b(9)\le157$ and $b(17)\le24337$. Further constructions
give, for example, $b(10)\le289$ and $b(13)\le2313$.
These bounds improve Lassak's bound $b(n)\le2^{n-1}+1$ of 1982 in every
dimension $n\ge4$. In dimension four they also improve the recent bounds
$b(4)\le8$ (Tolmachev and Voronov, May 2026) and $b(4)\le7$ (Arman,
Bondarenko, Prymak and Radchenko, October 2026); the bound $b(4)\le6$ was
announced independently by Chalmers in November 2026.

The bounds for $n\ge7$ have analytic proofs: a quantitative form of
Lassak's reduction to coverings of the sphere is combined with normal
partitions of weighted vector configurations, for which the smallest inner
product within a region is estimated explicitly. The bound for $n=6$ uses a
covering of $S^4$ by thirty sets, checked in rational arithmetic. For
$n=4,5$ three contact points of the circumscribed ball reduce the problem to
a four-parameter family of bodies, which is split by finitely many power
diagrams; that all parts have diameter less than one is verified by a
computer program with rigorous error control.
\end{abstract}

\medskip
\noindent\textbf{2020 Mathematics Subject Classification.}
Primary 52C17; Secondary 52A20, 65G30.

\noindent\textbf{Keywords.} Borsuk's problem, partition into sets of smaller
diameter, universal cover, contact points, covering of the sphere,
computer-assisted proof.

\section{Introduction}\label{sec:intro}

For a bounded set $A\subset\R^n$ of positive diameter let $b(A)$ be the
least number of sets of diameter strictly less than $\diam A$ whose union
is $A$, and let $b(n)$ be the supremum of $b(A)$ over all such $A$. Since
the overlaps of a covering can be removed one set at a time, it makes no
difference whether coverings or partitions are counted. We may also
normalize $\diam A=1$ and replace $A$ by its closure.

\subsection{Background}

In 1933 Borsuk~\cite{Borsuk} proved that a ball in $\R^n$ can be split into
$n+1$ parts of smaller diameter and asked whether $b(n)=n+1$ for all $n$.
The regular simplex shows that $b(n)\ge n+1$. The answer is positive for
$n\le3$; the three-dimensional case is due to Perkal (1947)~\cite{Perkal}
and Eggleston (1955)~\cite{Eggleston}, with simpler proofs by
Gr\"unbaum~\cite{Grunbaum} and Heppes~\cite{Heppes} (both 1957). In 1993
Kahn and Kalai~\cite{KahnKalai} showed that the conjecture fails in all
sufficiently large dimensions. The smallest dimension of a known
counterexample was later lowered to $298$ (2003)~\cite{HinrichsRichter},
$65$ (2014)~\cite{Bondarenko} and $64$ (2014)~\cite{JenrichBrouwer}. For
the history of the problem we refer to the surveys of
Raigorodskii~\cite{Raigorodskii,RaigorodskiiLectures},
Kalai~\cite{Kalai} and Zong~\cite{Zong}.

Two general upper bounds are relevant here. By Jung's
theorem (1901)~\cite{Jung}, every set of diameter one lies in a ball of radius
$\sqrt{n/(2n+2)}$. In 1982 Lassak~\cite{Lassak} intersected such a ball
with the unit ball centred at a point of contact and obtained
\begin{equation}\label{eq:lassak}
 b(n)\le2^{n-1}+1 .
\end{equation}
For large $n$, Schramm (1988)~\cite{Schramm} and, by a different method,
Bourgain and Lindenstrauss (1991)~\cite{BourgainLindenstrauss} proved
$b(n)\le(\sqrt{3/2}+o(1))^n$. For numerical comparison we use Schramm's
explicit estimate in the form given in~\cite[Theorem~6.2]{Naszodi},
\begin{equation}\label{eq:schramm}
 b(n)\le S(n):=5n^{3/2}(4+\log n)(3/2)^{n/2}.
\end{equation}
It is asymptotically far stronger than~\eqref{eq:lassak}, yet
$S(n)>2^{n-1}+1$ for $4\le n\le17$.

A closely related quantity is $d_{n,k}$, the supremum over all sets
$A\subset\R^n$ of diameter one of the infimum of those $d$ for which $A$
can be covered by $k$ sets of diameter at most $d$. Clearly $b(n)\le k$
whenever $d_{n,k}<1$. The planar numbers $d_{2,k}$ have been
studied since Lenz~\cite{Lenz} (1956); see Filimonov~\cite{Filimonov} and
Tolmachev, Protasov and Voronov~\cite{TPV}, who also showed that
$d_{3,4}\le0.966$. Bikeev and Raigorodskii~\cite{BikeevRaigorodskii}
recently studied the inverse question of how many parts of diameter at
most a given $b<1$ are needed. Our proofs give $d_{n,k}<1$ for every bound
$b(n)\le k$ obtained below (Remark~\ref{rem:dnk}).

In dimension four, Lassak's bound $b(4)\le9$ remained the best known for
more than forty years. In May 2026 Tolmachev and Voronov~\cite{TV} proved
$b(4)\le8$ by an optimized partition of truncated Lassak covers. In
October 2026 Arman, Bondarenko, Prymak and Radchenko~\cite{ABPR} proved
$b(4)\le7$ by a rigorous computer-assisted argument. They used contact
points of the circumscribed ball to reduce the problem to a four-parameter
family of bodies and covered every body of the family by seven balls of
diameter less than one. Our proof of $b(4)\le6$ uses a reduction of the
same type, and their facet-selection lemma (Lemma~\ref{lem:abpr} below) is
one of its ingredients, but the bodies are split in a different way: by
power diagrams into six cells, each of which is enclosed in a polytope of
diameter less than one. The same construction with twelve cells gives
$b(5)\le12$; the differences are listed at the beginning of
Section~\ref{sec:thm1}. A preprint of Chalmers~\cite{Chalmers}, posted on 4~November 2026,
also announces $b(4)\le6$. Our proof does not use that work, and we make
no claim of priority for the bound $b(4)\le6$.

\subsection{Main results}

\begin{mainthm}[dimensions four and five]\label{thm:main45}
$b(4)\le6$ and $b(5)\le12$.
\end{mainthm}

\begin{mainthm}[dimensions six to eight]\label{thm:main678}
$b(6)\le31$, $b(7)\le45$ and $b(8)\le89$.
\end{mainthm}

\begin{mainthm}[analytic estimate]\label{thm:analytic}
Let $n\ge9$. If $n-1=7a+8b+c$ with nonnegative integers $a,b,c$, then
\begin{equation}\label{eq:star-bound}
 b(n)\le1+88^a\,156^b\,2^c .
\end{equation}
In particular,
\begin{equation}\label{eq:eight}
 b(n)\le1+2^{n-1}\left(\tfrac{39}{64}\right)^{\lfloor(n-1)/8\rfloor}.
\end{equation}
Moreover, $b(10)\le289$.
\end{mainthm}

\begin{mainthm}[further bounds]\label{thm:ranges}
$b(12)\le1245$, $b(13)\le2313$ and $b(14)\le4977$.
\end{mainthm}

Theorems~\ref{thm:main45}--\ref{thm:analytic} improve Lassak's
bound~\eqref{eq:lassak} in every dimension $n\ge4$; for instance
$b(9)\le157$, whereas~\eqref{eq:lassak} gives $257$.
Theorem~\ref{thm:main45} is computer-assisted. The bound $b(6)\le31$ rests
on four $5\times5$ rational matrices and $63$ integer vectors, printed in
Appendix~\ref{app:rays}, which satisfy finitely many inequalities that are
checked in exact arithmetic. All other bounds have analytic proofs. The
general bound~\eqref{eq:eight} grows like $156^{n/8}\approx1.88^n$,
whereas Schramm's bound grows like $(3/2)^{n/2}$ up to a polynomial
factor, so our results concern small dimensions.

Table~\ref{tab:bounds} compares our estimates with the best previously
known ones. For $5\le n\le17$ this is Lassak's bound~\eqref{eq:lassak}.
In dimensions $18$ and $19$ Schramm's estimate~\eqref{eq:schramm} is
already stronger than~\eqref{eq:lassak}, and the table gives
$\lfloor S(n)\rfloor$; our bounds improve these values as well. In
dimension four the table gives the bound of~\cite{ABPR} and the bound
announced in~\cite{Chalmers}.

\begin{table}[ht]
\centering
\caption{Upper bounds for $b(n)$. In the column ``Previous'', the entries
marked with an asterisk~($^\ast$) are $\lfloor S(n)\rfloor$
from~\eqref{eq:schramm}, the entry marked with a dagger~($^\dagger$) is
the bound announced in~\cite{Chalmers}, the entry $7$ for $n=4$ is the
bound of~\cite{ABPR}, and the other entries are Lassak's bound
$2^{n-1}+1$. The last two columns
give the theorem that contains the bound and the type of its proof.}
\label{tab:bounds}
\begin{tabular}{@{}rrrll@{}}
\toprule
$n$ & Previous & This paper & Theorem & Type of proof\\\midrule
4&$7\;/\;6^\dagger$&6&\ref{thm:main45}&computer-assisted\\
5&17&12&\ref{thm:main45}&computer-assisted\\
6&33&31&\ref{thm:main678}&finite rational check\\
7&65&45&\ref{thm:main678}&analytic\\
8&129&89&\ref{thm:main678}&analytic\\
9&257&157&\ref{thm:analytic}&analytic\\
10&513&289&\ref{thm:analytic}&analytic\\
11&1025&625&\ref{thm:analytic}&analytic\\
12&2049&1245&\ref{thm:ranges}&analytic\\
13&4097&2313&\ref{thm:ranges}&analytic\\
14&8193&4977&\ref{thm:ranges}&analytic\\
15&16385&7745&\ref{thm:analytic}&analytic\\
16&32769&13729&\ref{thm:analytic}&analytic\\
17&65537&24337&\ref{thm:analytic}&analytic\\
18&$101144^\ast$&48673&\ref{thm:analytic}&analytic\\
19&$135395^\ast$&97345&\ref{thm:analytic}&analytic\\\bottomrule
\end{tabular}
\end{table}

\subsection{Methods}

The proofs use two different reductions. The first one, used for
Theorems~\ref{thm:main678}--\ref{thm:ranges}, refines Lassak's argument. After a set of diameter one is placed in his
lens, the lower part of the lens is cut into cones over the sets of a
covering of $S^{n-2}$. Lassak uses orthants, whose angular diameter is
$\pi/2$. We show that the angular diameter may exceed $\pi/2$ by an
explicit amount and determine the exact tolerance allowed by this argument
(Proposition~\ref{prop:lens}). It then remains to cover spheres by sets in
which all inner products are at least $-\mu$ for a small $\mu>0$. We use
the normal regions of weighted configurations that consist of half of the
vertices of a cube and some additional vectors: the coordinate vectors
(Proposition~\ref{prop:axis}), the vectors $2(\pm e_i\pm e_m)$ supported
on the last and one other coordinate (Proposition~\ref{prop:star}), or all
vectors supported on two coordinates (Proposition~\ref{prop:roots}). In
each case we bound the smallest inner product within a region. For the
second configuration all these inner products are nonnegative when $m=7$
or $m=8$, which gives coverings of $S^6$ by $88$ sets and of $S^7$ by $156$
sets, and hence Theorem~\ref{thm:analytic}. Coverings of
different dimensions are combined by an orthogonal product rule.

The second reduction, used for Theorem~\ref{thm:main45}, keeps three
contact points of the smallest enclosing ball. Following the idea
of~\cite{ABPR}, a contact triangle and a frame in the orthogonal complement
of its plane place every set of diameter one, up to an isometry, in a body
$H(u)$ from an explicit family with parameter $u\in[0,1]^4$; this works in
every dimension $n\ge4$ (Proposition~\ref{prop:normal3}). In dimensions four
and five we subdivide the parameter cube into dyadic boxes, replace the
bodies over each box by a rational outer approximation, and split it by a
rational power diagram into six or twelve parts. That every part has
diameter less than one is checked in exact arithmetic, or in floating point
with proved error bounds.

\subsection{Structure of the paper}\label{sec:structure}

Sections~\ref{sec:prelim} and~\ref{sec:spherical} contain the common
tools. Section~\ref{sec:prelim} fixes notation and recalls Jung's theorem.
Section~\ref{sec:spherical} reduces the problem to coverings of spheres
(Proposition~\ref{prop:lens}) and proves the lemmas used to bound inner
products within normal regions; it is not needed for
Theorem~\ref{thm:main45}. Each of Sections~\ref{sec:thm1}--\ref{sec:thm4}
proves one of Theorems~\ref{thm:main45}--\ref{thm:ranges}, and these four
sections can be read independently of one another (Table~\ref{tab:guide}).
Section~\ref{sec:table} explains how the entries of
Table~\ref{tab:bounds} are obtained from the coverings of spheres
constructed in Sections~\ref{sec:thm2}--\ref{sec:thm4}, and
Section~\ref{sec:remarks} contains concluding remarks.

The appendices describe data and computations. Appendix~\ref{app:rays}
lists the rational data used in Section~\ref{sec:six}.
Appendix~\ref{app:verify} explains how the partitions of
Section~\ref{sec:thm1} are verified, and Appendix~\ref{app:search} how
they were found; the search plays no role in the proof. The certificates and verification programs
are available in the repository~\cite{Repo}.

\begin{table}[ht]
\centering
\caption{Where the main results are proved.}
\label{tab:guide}
\begin{tabular}{@{}lllll@{}}
\toprule
Result & Bounds & Proof & Uses & Computations\\\midrule
Theorem~\ref{thm:main45} & $n=4,5$ & Section~\ref{sec:trees} & Section~\ref{sec:prelim} & Appendices~\ref{app:verify}, \ref{app:search}\\
Theorem~\ref{thm:main678} & $n=6,7,8$ & Section~\ref{sec:proof2} & Sections~\ref{sec:prelim}, \ref{sec:spherical} & Appendix~\ref{app:rays} ($n=6$)\\
Theorem~\ref{thm:analytic} & $n\ge9$ & Section~\ref{sec:proof3} & Sections~\ref{sec:prelim}, \ref{sec:spherical} & none\\
Theorem~\ref{thm:ranges} & $n=12,13,14$ & Section~\ref{sec:proof4} & Sections~\ref{sec:prelim}, \ref{sec:spherical} & none\\\bottomrule
\end{tabular}
\end{table}

\section{Notation and common facts}\label{sec:prelim}

We write $\ip xy$ and $|x|$ for
the Euclidean inner product and norm, $e_1,\dots,e_n$ for the standard
basis, $B(c,r)$ for the closed ball with centre $c$ and radius $r$, and
$S^{m-1}$ for the unit sphere of $\R^m$. For $C\subset\R^m$ we put
$\cone C=\{tx:t\ge0,\ x\in C\}$. The support function of a bounded set $M$
is $h_M(v)=\sup_{x\in M}\ip vx$. We set $\diam\varnothing=0$. If $\diam M\le1$, then for every unit
vector $v$
\begin{equation}\label{eq:width}
 h_M(v)+h_M(-v)=\sup_{x,y\in M}\ip v{x-y}\le1.
\end{equation}

Every compact set $X\subset\R^n$ with more than one point has a unique
smallest enclosing ball, the \emph{circumscribed ball} of $X$; its radius
is the \emph{circumradius} of $X$~\cite{DGK}. Throughout the paper
\[
 r_n=\sqrt{\frac n{2n+2}},\qquad \rho_n=r_n^2=\frac n{2n+2}.
\]

\begin{theorem}[Jung, 1901~\cite{Jung}]\label{thm:jung}
If $X\subset\R^n$ is compact and $\diam X=1$, then the circumradius of $X$
lies between $1/2$ and $r_n$.
\end{theorem}

We also use Borsuk's theorem~\cite{Borsuk} that a ball in $\R^n$ can be
partitioned into $n+1$ sets of smaller diameter.

\section{Coverings of spheres and Lassak's lens}\label{sec:spherical}

This section collects the general tools for the proofs of
Theorems~\ref{thm:main678}, \ref{thm:analytic} and~\ref{thm:ranges}; the
proof of Theorem~\ref{thm:main45} does not use it. The main statement is
Proposition~\ref{prop:lens}, which reduces upper bounds for $b(n)$ to
coverings of the sphere $S^{n-2}$; Corollary~\ref{cor:blocks} combines
coverings of spheres of different dimensions.

\subsection{Coverings of spheres}\label{sec:coverings}

\begin{definition}\label{def:K}
For $m\ge1$ and $0\le\mu<1$ let $K(m,\mu)$ be the least number of closed
sets covering $S^{m-1}$ such that $\ip uv\ge-\mu$ whenever $u$ and $v$
belong to the same set.
\end{definition}

Clearly $K(m,\mu)$ is nonincreasing in $\mu$. The closed orthants give
$K(m,0)\le2^m$, and $K(1,0)=2$.

\begin{lemma}\label{lem:product}
For $a,b\ge1$ and $0\le\mu<1$ we have $K(a+b,\mu)\le K(a,\mu)K(b,\mu)$.
\end{lemma}

\begin{proof}
Let $C_1,\dots,C_M$ and $D_1,\dots,D_N$ be coverings of $S^{a-1}$ and
$S^{b-1}$ as in Definition~\ref{def:K}. The sets
\[
 U_{ij}=S^{a+b-1}\cap(\cone C_i\times\cone D_j)
\]
are closed and cover $S^{a+b-1}$. For two points $(u,v)$ and $(u',v')$ of
the same set $U_{ij}$, the Cauchy--Schwarz inequality gives
\[
 \ip{(u,v)}{(u',v')}\ge-\mu\bigl(|u||u'|+|v||v'|\bigr)\ge-\mu. \qedhere
\]
\end{proof}

\subsection{Lassak's lens and products of coverings}\label{sec:lens}

For $n\ge2$ we define
\begin{equation}\label{eq:muhat}
 s_n=\sqrt{\frac{n-1}{2n}},\qquad
 \muhat_n=\frac{1-s_n}{(n+1)s_n}
 =\frac{\sqrt{2n/(n-1)}-1}{n+1}.
\end{equation}

\begin{proposition}[quantitative Lassak reduction]\label{prop:lens}
Let $n\ge2$ and $0\le\mu<\muhat_n$. Then
\[
 b(n)\le1+K(n-1,\mu).
\]
\end{proposition}

\begin{figure}[t]
\centering
\includegraphics[width=\textwidth]{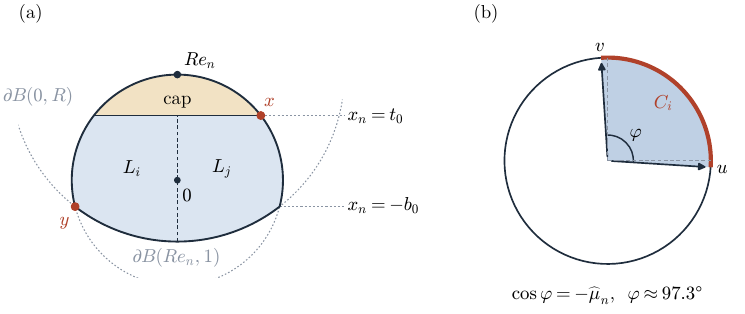}
\caption{(a) Meridian section of Lassak's lens $L=B(0,R)\cap B(Re_n,1)$,
$R=r_n$, for $n=4$. The part of $L$ above the height $t_0$ (the cap) lies in
a ball of radius $\sqrt{R^2-t_0^2}=1/2$; cutting slightly higher, at
$t>t_0$, makes its diameter less than one. Below the cut, $L$ is divided
into the pieces $L_i$ over the sets $C_i$ of a covering of $S^{n-2}$. The
extremal points $x=(u/2,t_0)$ and $y=(s_nv,-b_0)$ of Remark~\ref{rem:lens}
lie on $\partial L$, and $y$ lies on the rim where the two spheres meet;
the section shows only their heights and distances from the axis.
(b) The angular diameter of a set $C_i$ may be any number less than
$\varphi$, where $\cos\varphi=-\widehat\mu_n$; this is slightly more than
the right angle of an orthant (dashed).}
\label{fig:lens}
\end{figure}

For $\mu=0$, together with $K(n-1,0)\le2^{n-1}$, this is Lassak's
bound~\eqref{eq:lassak}. The point of Proposition~\ref{prop:lens} is that the
sets of the spherical covering may have angular diameter slightly larger
than $\pi/2$ (Figure~\ref{fig:lens}(b)).

\begin{proof}
Let $A\subset\R^n$ be compact with $\diam A=1$, and write $R=r_n$ for
brevity. By Theorem~\ref{thm:jung}, after a translation the circumscribed
ball of $A$ is $B(0,r)$ with $r\le R$. Choose a point of $A$ on its
boundary and rotate it to $re_n$. The ball $B((r-R)e_n,R)$ contains
$B(0,r)$ and touches it at $re_n$. Translating its centre to the origin, we
obtain
\begin{equation}\label{eq:lens}
 A\subset L:=B(0,R)\cap B(Re_n,1),\qquad Re_n\in A,
\end{equation}
where the second ball comes from $\diam A=1$. This is Lassak's lens; see
Figure~\ref{fig:lens}(a). Put
\[
 \delta=1-2R^2=\frac1{n+1},\qquad s=s_n,\qquad
 t_0=Rs=\sqrt{R^2-\tfrac14},\qquad b_0=\frac{\delta}{2R}.
\]
A direct computation gives
\begin{equation}\label{eq:lens-identities}
 R^2-b_0^2=s^2,\qquad
 \muhat_n=\frac{\delta(1-s)}s<\frac1n,\qquad
 \frac{b_0}{2t_0s}=\frac1{n-1},\qquad
 \frac{2t_0s}{b_0}=n-1.
\end{equation}
The strict inequality is equivalent to $s>n/(2n+1)$, which holds because
$s\ge1/2$ for $n\ge2$.

Let $C_1,\dots,C_N$, $N=K(n-1,\mu)$, be closed sets that cover $S^{n-2}$
and in which all inner products are at least $-\mu$. We write points of
$\R^n$ as $(x',x_n)$ with $x'\in\R^{n-1}$.

\emph{Step 1: the lower part below height $t_0$.} Let $x=(\xi u,\alpha)$ and
$y=(\zeta v,\beta)$, where $\xi,\zeta\ge0$, $u,v\in C_i$ for the same $i$,
and $x,y\in L\cap\{x_n\le t_0\}$. (If $\xi=0$ or $\zeta=0$, the direction
$u$ or $v$ may be chosen arbitrarily.) Since $\ip uv\ge-\mu$,
\begin{equation}\label{eq:lens-distance}
 |x-y|^2\le |x|^2+|y|^2-2\alpha\beta+2\mu\xi\zeta.
\end{equation}
If $\alpha\beta\ge0$, the right-hand side is at most $2R^2(1+\mu)$, which
is less than one because $\mu<1/n$. Otherwise we may assume that
$0<\alpha\le t_0$ and $\beta=-b<0$.

Suppose first that $0<b\le b_0$. Using $|x|,|y|\le R$,
$\xi\le\sqrt{R^2-\alpha^2}$ and $\zeta\le\sqrt{R^2-b^2}$, we bound the
right-hand side of~\eqref{eq:lens-distance} by
\[
 F(\alpha,b)=2R^2+2\alpha b+
 2\mu\sqrt{R^2-\alpha^2}\sqrt{R^2-b^2},
 \qquad (\alpha,b)\in[0,t_0]\times[0,b_0].
\]
For $0<\mu<1$ the function $F$ has no critical point with $\alpha,b>0$:
the equations $F_\alpha=F_b=0$ multiply to $\mu^2=1$. For $\mu=0$ it is
bilinear. Hence $F$ attains its maximum on the boundary of the rectangle.
On the sides $\alpha=0$ and $b=0$ we have $F\le2R^2(1+\mu)$. On the other
two sides, by~\eqref{eq:lens-identities} and $\sqrt{R^2-t_0^2}=1/2$,
\[
 \tfrac12 F_\alpha(\alpha,b_0)\ge b_0-2\mu st_0>0,
 \qquad
 F_b(t_0,b)\ge2t_0-\mu b_0/s>0.
\]
Therefore the maximum over these two sides is $F(t_0,b_0)=2R^2+\delta
s+\mu s$, and
\begin{equation}\label{eq:lens-M}
 F(\alpha,b)\le
 M:=\max\{2R^2(1+\mu),\ 2R^2+\delta s+\mu s\}<1.
\end{equation}
The second term in~\eqref{eq:lens-M} is less than one exactly when
$\mu<\delta(1-s)/s=\muhat_n$.

Now let $b\ge b_0$. The ball $B(Re_n,1)$ gives $|y|^2\le1-R^2-2Rb$ and
$\zeta\le\sqrt{1-(R+b)^2}$, and therefore
\[
 |x-y|^2\le1-2(R-\alpha)b+
 2\mu\sqrt{R^2-\alpha^2}\sqrt{1-(R+b)^2}.
\]
For fixed $\alpha\le t_0<R$ the right-hand side decreases in $b$. Since
$2Rb_0=\delta=1-2R^2$, its value at $b=b_0$ equals $F(\alpha,b_0)\le M$.

\emph{Step 2: thickening the lower parts.} Put
\begin{equation}\label{eq:lens-eps}
 \eps=\min\left\{\frac{R-t_0}{2},\ \frac{1-\sqrt M}{4}\right\}>0,
 \qquad t=t_0+\eps,
\end{equation}
and $L_i(t)=\{x\in L:x_n\le t,\ x'\in\cone C_i\}$. The map
$T(x',x_n)=(x',\min\{x_n,t_0\})$ moves each point by at most $\eps$ and
sends $L_i(t)$ into $L_i(t_0)$. Indeed, if $x_n>t_0$, then $|Tx|\le|x|$
and $|x'|^2\le R^2-t_0^2$, so $|Tx-Re_n|^2\le2R^2-2Rt_0<1$. Consequently
\[
 \diam L_i(t)\le\sqrt M+2\eps\le\frac{1+\sqrt M}{2}<1.
\]

\emph{Step 3: the cap.} The remaining part $L\cap\{x_n\ge t\}$ lies in the
ball $B\bigl(te_n,\sqrt{R^2-t^2}\bigr)$, so its diameter is at most
$2\sqrt{R^2-t^2}$, which is less than one because $t>t_0$.

The $N+1$ sets $L_1(t),\dots,L_N(t)$ and $L\cap\{x_n\ge t\}$ cover $A$ and
have diameter less than one.
\end{proof}

\begin{remark}\label{rem:lens}
The proof gives a uniform bound: every part has diameter at most
\begin{equation}\label{eq:uniform-diameter}
 \Delta(n,\mu)=\max\left\{\frac{1+\sqrt M}{2},\
 2\sqrt{r_n^2-(t_0+\eps)^2}\right\}<1,
\end{equation}
which depends only on $n$ and $\mu$. The threshold $\muhat_n$ cannot be
improved within this argument. For $n\ge3$ take $u,v\in S^{n-2}$ with
$\ip uv=-\muhat_n$. The points $x=(u/2,t_0)$ and $y=(s_nv,-b_0)$ lie in
$L$, and $|x-y|^2=2r_n^2+\delta s_n+\muhat_ns_n=1$; both points remain in
the lower part if the truncation height is raised. This says nothing, of
course, about the sharpness of the resulting bounds for $b(n)$.
\end{remark}

To check the hypothesis of Proposition~\ref{prop:lens} for a rational $\mu$ we
use the following equivalent form: if $\mu=p/q$ with integers $p\ge0$ and
$q>0$, then
\begin{equation}\label{eq:slack}
 \mu<\muhat_n\quad\Longleftrightarrow\quad
 D(n;p,q):=2nq^2-(n-1)\bigl(q+(n+1)p\bigr)^2>0.
\end{equation}
Indeed, $\mu<\muhat_n$ means $s_n\bigl(1+(n+1)\mu\bigr)<1$; it remains to
square both sides. Note also that $\muhat_n$ is strictly decreasing in $n$.

We call a triple $(m,N,\mu)$ with $K(m,\mu)\le N$ a \emph{block}. The
following corollary combines blocks of different dimensions.

\begin{corollary}\label{cor:blocks}
Let $(m_j,N_j,\mu_j)$ be a finite list of blocks that includes $(1,2,0)$.
For $n\ge2$ let $I_n=\{j:\mu_j<\muhat_n\}$. Then
\begin{equation}\label{eq:block-formula}
 b(n)\le1+
 \min\Bigl\{\prod_{j\in I_n}N_j^{q_j}:\ q_j\in\Z_{\ge0},\
 \sum_{j\in I_n}m_jq_j=n-1\Bigr\}.
\end{equation}
\end{corollary}

\begin{proof}
For an admissible choice of the $q_j$ apply Lemma~\ref{lem:product} with
$\mu=\max\{\mu_j:q_j>0\}$, using that $K(m,\mu)$ is nonincreasing in
$\mu$, and then Proposition~\ref{prop:lens}. The one-dimensional block
guarantees that the minimum is taken over a nonempty set.
\end{proof}

The minimum in~\eqref{eq:block-formula} is over finitely many integer
vectors and is computed by elementary dynamic programming.

\subsection{Normal regions and four elementary lemmas}\label{sec:elementary}

\begin{figure}[t]
\centering
\includegraphics[width=0.42\textwidth]{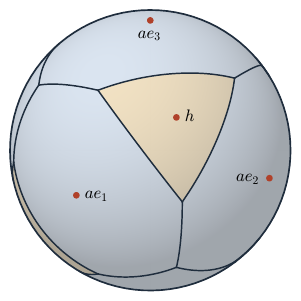}
\caption{Normal regions of $\{\pm2e_i\}\cup H_3$ on $S^2$, the
three-dimensional analogue of the configurations $V^A_{m,a}$ of
Section~\ref{sec:axis} (where $m\ge6$). The regions of the coordinate
vectors are light, those of the half-cube are shaded; the boundaries are
the radial projections of the edges of the polar polytope of the
configuration.}
\label{fig:regions}
\end{figure}

All our coverings come from finite configurations of vectors. Let
$V\subset\R^m$ be a finite set that spans $\R^m$. The \emph{normal
regions} of $V$ are
\begin{equation}\label{eq:normal-region}
 C_v(V)=\{x\in S^{m-1}:\ip xv\ge\ip xw\text{ for every }w\in V\}
 \qquad(v\in V).
\end{equation}
They are closed and cover the sphere, so $K(m,\mu)\le|V|$ as soon as all
inner products within each region are at least $-\mu$. The lengths of the
vectors matter: rescaling some of them changes the regions. When the
origin is an interior point of $\conv V$, the cones over the regions of
the vertices $v$ of $\conv V$ are the normal cones of $\conv V$ at these
vertices, that is, the cones over the facets of the polar polytope. Figure~\ref{fig:regions} shows
a three-dimensional example.

The configurations in Sections~\ref{sec:axis}, \ref{sec:thm3}
and~\ref{sec:thm4} contain the even half of the vertex set of the cube,
\[
 H_m=\Bigl\{\varepsilon\in\{-1,1\}^m:\ \prod_i\varepsilon_i=1\Bigr\},
 \qquad h=(1,\dots,1)\in H_m .
\]
For $x\in\R^m$ we write $S_x=\sum_ix_i$.

The smallest inner products in the normal regions of these
configurations are computed with the following four lemmas. The first is used to maximize convex quadratic forms.

\begin{lemma}\label{lem:vertices}
Let $w_1,\dots,w_N\ge0$ and $q>0$, and let
\[
 P=\Bigl\{z\in\R^N:\ z_i\ge0,\ \sum_iz_i\le1,\ \sum_iw_iz_i\le q\Bigr\}.
\]
Every convex function on $P$ attains its maximum at a vertex of $P$. Every
vertex $z\ne0$ of $P$ has one of the following two forms:
\begin{enumerate}[label=\textup{(\alph*)}]
\item exactly one coordinate $z_i$ is nonzero, and $z_i=\min\{1,q/w_i\}$
(with $q/0=\infty$);
\item exactly two coordinates $z_i,z_j$ are nonzero, $w_i\ne w_j$, and
$\sum_iz_i=1$, $\sum_iw_iz_i=q$.
\end{enumerate}
\end{lemma}

\begin{proof}
Since $0\le z_i\le1$ on $P$, the set $P$ is a polytope, and a convex
function attains its maximum over a polytope at a
vertex~\cite[Section~32]{Rockafellar}. At a vertex, $N$ linearly
independent constraints are active. At most two of them are not of the
form $z_i\ge0$, so at most two coordinates are nonzero. If $z_i$ and $z_j$
are nonzero, both remaining constraints must be active, and their
restrictions $(1,1)$ and $(w_i,w_j)$ to these coordinates must be linearly
independent; hence $w_i\ne w_j$. If only $z_i$ is nonzero, one of the two
constraints is active and has a nonzero coefficient at $z_i$, so
$z_i\in\{1,q/w_i\}$; feasibility gives $z_i=\min\{1,q/w_i\}$.
\end{proof}

\begin{lemma}[caps]\label{lem:caps}
Let $w\in S^{m-1}$, $0<\gamma\le1$, and let $C\subset S^{m-1}$ satisfy
$\ip xw\ge\sqrt\gamma$ for all $x\in C$. Then $\ip xy\ge2\gamma-1$ for
all $x,y\in C$.
\end{lemma}

\begin{proof}
Let $\theta=\arccos\sqrt\gamma\in[0,\pi/2)$. The angles between $w$ and
$x$, and between $w$ and $y$, are at most $\theta$. By the triangle
inequality for the angular metric, the angle between $x$ and $y$ is at most
$2\theta<\pi$, so $\ip xy\ge\cos2\theta=2\gamma-1$.
\end{proof}

\begin{lemma}[sign cones]\label{lem:signcone}
Let $x,y\in\R^m$ satisfy $x_i+x_j\ge0$ and $y_i+y_j\ge0$ for all $i\ne j$,
and $2|x_i|\le S_x$, $2|y_i|\le S_y$ for all $i$. Then each of $x$ and $y$
has at most one negative coordinate. If $\ip xy<0$, then both $x$ and $y$
have a negative coordinate, and these coordinates are in different
positions.
\end{lemma}

\begin{proof}
Two negative coordinates would violate $x_i+x_j\ge0$. If neither vector
has a negative coordinate, then $\ip xy\ge0$. If both have their negative
coordinate in the same position $i$, then $x_iy_i>0$ and all other terms of
$\ip xy$ are nonnegative. Suppose that only $x$ has a negative coordinate,
say $x_i=-t<0$. Then $x_j\ge t$ for $j\ne i$, all $y_j$ are nonnegative,
and
\[
 \ip xy=-ty_i+\sum_{j\ne i}x_jy_j\ge t\Bigl(\sum_{j\ne i}y_j-y_i\Bigr)
 =t(S_y-2y_i)\ge0 .
\]
The case where only $y$ has a negative coordinate is symmetric.
\end{proof}

\begin{lemma}[negative pairs]\label{lem:negative-pair}
Let $k\ge1$, $c\ge1$ and $1\le\ell<k$. Suppose that
\[
 x=(-1,b,1+p_1,\ldots,1+p_k),\qquad
 y=(d,-1,1+q_1,\ldots,1+q_k),
\]
where $b,d\ge1$, $p_i,q_i\ge0$, $\sum p_i\ge c(b-\ell)_+$ and
$\sum q_i\ge c(d-\ell)_+$. Then
\[
 \frac{\ip xy}{|x||y|}\ge
 -\max\left\{0,\frac{2\ell-k}{\ell^2+k+1}\right\}.
\]
\end{lemma}

Since $p_i,q_i\ge0$, the hypotheses on $\sum p_i$ and $\sum q_i$ are
equivalent to $\sum p_i\ge c(b-\ell)$ and $\sum q_i\ge c(d-\ell)$.

\begin{proof}
Let $f(z)=z-c(z-\ell)_+$; since $c\ge1$, we have $f(z)\le\min\{z,\ell\}$.
Hence
\[
 \ip xy\ge k-b-d+\sum p_i+\sum q_i\ge k-f(b)-f(d)\ge k-b'-d',
\]
where $b'=\min\{b,\ell\}$ and $d'=\min\{d,\ell\}$. If $k-b'-d'\ge0$ there
is nothing to prove. Otherwise the bounds $|x|^2\ge b'^2+k+1$ and
$|y|^2\ge d'^2+k+1$ give
\[
 \frac{\ip xy}{|x||y|}\ge
 \frac{k-b'-d'}{\sqrt{b'^2+k+1}\,\sqrt{d'^2+k+1}}.
\]
The right-hand side decreases in each variable on $[1,\ell]$: its partial
derivative in $b'$ has the sign of $-(k+1)-b'(k-d')<0$. The minimum is
therefore attained at $b'=d'=\ell$.
\end{proof}

\section{Proof of Theorem~\ref{thm:main45}: three contact points}\label{sec:thm1}

Following the idea
of~\cite{ABPR}, we fix more of the geometry of a set of diameter one than
Lassak's lens does. The proof has three steps.
\begin{enumerate}[label=\arabic*.]
\item Three contact points of the circumscribed ball and a suitable frame
place every set of diameter one, up to an isometry, in a body $H(u)$ from
an explicit family with parameter $u\in[0,1]^4$
(Proposition~\ref{prop:normal3}).
\item For each box $Q$ of a dyadic subdivision of $[0,1]^4$, the bodies
$H(u)$ with $u\in Q$ lie in one body $\Omega_Q$ with rational data, so it
suffices to split finitely many bodies $\Omega_Q$
(Proposition~\ref{prop:reduction3}).
\item Each $\Omega_Q$ is split by a rational power diagram into six parts
($n=4$) or twelve parts ($n=5$), and every part is enclosed in a polytope of
diameter less than one (Definition~\ref{def:admissible} and
Proposition~\ref{prop:admissible}).
\end{enumerate}
The verification evaluates the enclosures of Step~2 in
outward-rounded interval arithmetic and checks, for every box, that the
polytopes of Step~3 have diameter less than one. The checking program and
its error analysis are described in Appendix~\ref{app:verify}.

\paragraph{Relation to~\cite{ABPR}.} Three ingredients are taken
from~\cite{ABPR}: the idea of reducing the problem, by means of contact
points of the circumscribed ball, to a family of bodies with four
parameters; the facet-selection lemma (Lemma~\ref{lem:abpr}); and the
interval arithmetic used to evaluate the enclosures
(Section~\ref{sec:boxes}). The splitting of the bodies is different.
In~\cite{ABPR} every body of the family is covered by seven balls of
diameter less than one. Here the cells are those of a power diagram: they
are polytopes, adapted to the shape of $\Omega_Q$, and their number is six
for $n=4$. The parametrization by a contact triangle and a frame
(Definition~\ref{def:H}) is valid in every dimension $n\ge4$, and we use
it also for $n=5$. The partition trees, the conditions of
Definition~\ref{def:admissible} and the programs that check them are our
own (Appendices~\ref{app:verify} and~\ref{app:search}).

\subsection{Two auxiliary lemmas}\label{sec:contact-lemmas}

We first record the contact-point criterion for a smallest enclosing ball.

\begin{lemma}\label{lem:center}
Let $B(c,r)$ be the circumscribed ball of a compact set $X\subset\R^n$,
and let $T=X\cap\partial B(c,r)$. Then $c\in\conv T$.
\end{lemma}

\begin{proof}
The set $T$ is compact, and it is nonempty since otherwise the ball could
be shrunk. Suppose that $c\notin\conv T$. Then $c$ is strictly separated
from the compact convex set $\conv T$: there are a unit vector $u$ and
$\delta>0$ with $\ip u{t-c}\ge\delta$ for all $t\in T$. Let
$U=\{x\in X:\ip u{x-c}>\delta/2\}$. The compact set $X\setminus U$ does not
meet $T$, so $r'=\max_{x\in X\setminus U}|x-c|<r$ (with $r'=0$ if
$X\setminus U=\varnothing$). Let $c_\eps=c+\eps u$ with
$0<\eps<\min\{\delta/2,\,r-r'\}$. For $x\in U$,
\[
 |x-c_\eps|^2=|x-c|^2-2\eps\ip u{x-c}+\eps^2< r^2-\eps\delta+\eps^2<r^2,
\]
and for $x\in X\setminus U$ we have $|x-c_\eps|\le r'+\eps<r$. Thus $X$
lies in a ball of radius less than $r$, a contradiction.
\end{proof}

\begin{lemma}[frame]\label{lem:frame}
Let $W$ be a Euclidean space of dimension $d\ge1$, and let $M\subset W$ be
compact with $\diam M\le1$. Then $W$ has an orthonormal basis
$v_1,\dots,v_d$ such that
\[
 h_M(v_i)\le\tfrac12,\quad h_M(-v_i)\le\tfrac12\quad (1\le i\le d-1),
 \qquad h_M(v_d)\le\tfrac12 .
\]
\end{lemma}

\begin{proof}
The function $f(v)=h_M(v)-h_M(-v)$ is continuous (support functions of
compact sets are Lipschitz) and odd. Suppose that $v_1,\dots,v_{j-1}$ have
been chosen, $j\le d-1$. Their orthogonal complement has dimension
$d-j+1\ge2$. On its unit sphere, $f$ vanishes somewhere: if $f(v)\ne0$,
apply the intermediate value theorem on a great semicircle from $v$ to
$-v$. This gives a unit vector $v_j$ with $h_M(v_j)=h_M(-v_j)$; by~\eqref{eq:width} both values are at
most $1/2$. Finally, $v_d$ is determined up to sign, and we choose the sign
with $h_M(v_d)\le h_M(-v_d)$; then $h_M(v_d)\le1/2$, again
by~\eqref{eq:width}.
\end{proof}

\subsection{The contact triangle}\label{sec:triangle}

Put
\[
 c_n=\prod_{d=3}^{n}\Bigl(1-\frac1{d^2}\Bigr)=\frac{2(n+1)}{3n},
 \qquad\text{so that } c_4=\tfrac56,\quad c_5=\tfrac45,\quad c_n\rho_n=\tfrac13.
\]
The product telescopes, because
$\prod_{d=2}^{n}(1-1/d^2)=\prod_{d=2}^{n}\frac{(d-1)(d+1)}{d^2}=\frac{n+1}{2n}$.
We use the following lemma of Arman, Bondarenko, Prymak and Radchenko.

\begin{lemma}[{\cite[Lemma~3]{ABPR}}]\label{lem:abpr}
Let $d\ge2$, and let $\Delta$ be a $d$-simplex whose circumcentre lies in
its relative interior. Then $\Delta$ has a facet whose circumcentre lies in
the relative interior of the facet and whose squared circumradius is at
least $1-1/d^2$ times the squared circumradius of $\Delta$.
\end{lemma}

\begin{proof}
Work in the affine hull of $\Delta$, put its circumcentre at $0$, and
write its vertices as $v_0,\ldots,v_d$, with $|v_i|=R$.
Let $\theta_i>0$ be the barycentric coordinates of $0$.
Write the facet opposite $v_i$ as $\ip{n_i}x=h_i$, where $|n_i|=1$,
$h_i>0$ and $\Delta\subset\{x:\ip{n_i}x\le h_i\}$.
Choose $i$ with $\theta_i\le1/(d+1)$. Taking the scalar product of
$\sum_j\theta_jv_j=0$ with $n_i$ gives
\[
 h_i=-\frac{\theta_i}{1-\theta_i}\ip{n_i}{v_i}\le\frac Rd.
\]
Now choose $j$ for which $h_j$ is smallest. The point $h_jn_j$ lies in
the relative interior of its facet: for every $l\ne j$,
$\ip{n_l}{h_jn_j}<h_j\le h_l$, since the distinct facet normals are unit
vectors. This point is the orthogonal projection of $0$ onto that facet,
and hence its circumcentre. The squared facet circumradius is
$R^2-h_j^2\ge(1-1/d^2)R^2$.
\end{proof}

\begin{lemma}[contact triangle]\label{lem:triangle}
Let $X\subset\R^n$ be compact with $\diam X=1$, let $B(0,\sqrt\rho)$ be
its circumscribed ball, and suppose that $\rho>1/4$. Then there are points $P_1,P_2,P_3\in X$ with $|P_i|^2=\rho$
that form an acute triangle, and its squared circumradius $\tau$ satisfies
\[
 c_n\rho\le\tau\le\min\{\rho,\,1/3\}.
\]
Moreover, the circumcentre of the triangle is the orthogonal projection of
$0$ onto its plane and lies at distance $h=\sqrt{\rho-\tau}$ from $0$.
\end{lemma}

\begin{proof}
Let $T=X\cap\partial B(0,\sqrt\rho)$. By Lemma~\ref{lem:center} and
Carath\'eodory's theorem there is a finite set $T_0\subset T$ with
$0\in\conv T_0$; take one of smallest cardinality. By minimality, the
points of $T_0$ are affinely independent and all barycentric coordinates
of $0$ with respect to them are positive. Hence $\Delta_0=\conv T_0$ is a
$d_0$-simplex, $d_0\le n$, with $0$ in its relative interior. Its vertices
lie on the sphere $|x|^2=\rho$, whose centre lies in $\aff\Delta_0$, so $0$
is the circumcentre of $\Delta_0$ and its squared circumradius is $\rho$.
Here $d_0\ge2$: the case $d_0=0$ is impossible because $\rho>0$, and
$d_0=1$ would mean $T_0=\{p,-p\}$ with $|2p|=2\sqrt\rho>1$.

For a face $F$ of $\Delta_0$, the point of $\aff F$ nearest to $0$ is
equidistant from the vertices of $F$. It is therefore the circumcentre of
$F$, and the squared circumradius of $F$ equals $\rho$ minus the squared
distance from $0$ to $\aff F$. Apply Lemma~\ref{lem:abpr} to $\Delta_0$
with $d=d_0$, then to the facet obtained, regarded as a simplex in its own
affine hull, with $d=d_0-1$, and so on down to $d=3$. The result is a
triangle $P_1P_2P_3$ with vertices in $T$ whose circumcentre lies in its
relative interior and whose squared circumradius satisfies
\[
 \tau\ge\prod_{d=3}^{d_0}\Bigl(1-\frac1{d^2}\Bigr)\rho\ge c_n\rho .
\]
By the previous paragraph $\tau=\rho-h^2\le\rho$. The triangle is acute,
so its largest angle $\alpha$ satisfies $60^\circ\le \alpha<90^\circ$, and
its longest side equals $2\sqrt\tau\sin \alpha\ge\sqrt{3\tau}$. Since all
sides are at most one, $\tau\le1/3$.
\end{proof}

The two bounds on $\tau$ are compatible exactly up to Jung's bound:
$c_n\rho\le1/3$ if and only if $\rho\le\rho_n$.

\subsection{The normal form}\label{sec:normalform}

\begin{definition}\label{def:H}
For $u=(u_1,u_2,u_3,u_4)\in[0,1]^4$ put
\begin{align*}
t&=u_1\sqrt{\rho_n-1/3}-(1-u_1)\sqrt{1/12}, & \rho&=1/3+t|t|,\\
h&=(1-u_2)\sqrt{(1-c_n)\rho}+u_2\max\{t,0\}, & \tau&=\rho-h^2,\\
a_0&=\sqrt{\max\{0,\,1-1/(4\tau)\}}, & a&=(1-u_3)/2+u_3a_0,\\
b&=(1-u_4)\sqrt{(1-a)/2}+u_4a, &&
\end{align*}
and
\begin{align*}
P_1&=\bigl(\sqrt\tau,\ 0,\ h,\ 0,\dots,0\bigr),\\
P_2&=\bigl(\sqrt\tau\,(2a^2-1),\ 2\sqrt\tau\,a\sqrt{1-a^2},\ h,\ 0,\dots,0\bigr),\\
P_3&=\bigl(\sqrt\tau\,(2b^2-1),\ -2\sqrt\tau\,b\sqrt{1-b^2},\ h,\ 0,\dots,0\bigr).
\end{align*}
Let
\[
 F_n=\{x\in\R^n:\ |x_j|\le\tfrac12\ (4\le j\le n-1),\ x_n\le\tfrac12\},\qquad
 H(u)=B(0,\sqrt\rho)\cap\bigcap_{i=1}^3B(P_i,1)\cap F_n .
\]
\end{definition}

\begin{figure}[t]
\centering
\includegraphics[width=\textwidth]{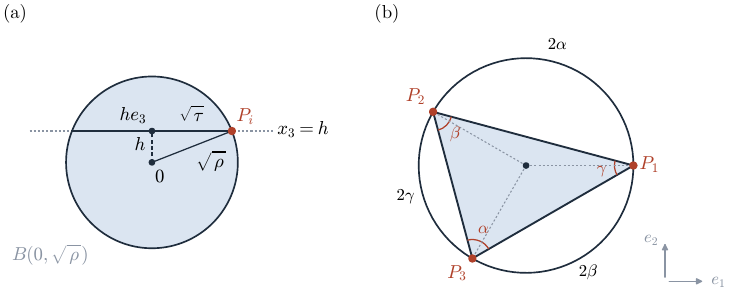}
\caption{The normal form of Definition~\ref{def:H}, drawn for $\rho=0.3$,
$\tau=0.26$, $\alpha=75^\circ$, $\beta=60^\circ$. (a) Section through
$e_3$: the plane of the contact triangle is at height $h=\sqrt{\rho-\tau}$
and meets the sphere $|x|^2=\rho$ in a circle of radius $\sqrt\tau$.
(b) The plane $x_3=h$ seen from above: $P_2$ and $P_3$ are obtained from
$P_1$ by rotations through $2\alpha$ and $-2\beta$, where
$\alpha\ge\beta\ge\gamma$ are the angles of the triangle at $P_3$, $P_2$,
$P_1$.}
\label{fig:triangle}
\end{figure}

Figure~\ref{fig:triangle} illustrates this parameterization. To check that
it is well defined, note that as $u_1$ runs over $[0,1]$, $t$ increases from
$-\sqrt{1/12}$ to $\sqrt{\rho_n-1/3}$, and $\rho$ increases continuously
from $1/4$ to $\rho_n$. Next,
$\max\{t,0\}=\sqrt{\max\{0,\rho-1/3\}}\le\sqrt{(1-c_n)\rho}$, because
$\rho-1/3\le(1-c_n)\rho$ is equivalent to $c_n\rho\le1/3$. Hence $h$ is a
convex combination of two numbers in $[\max\{t,0\},\sqrt{(1-c_n)\rho}]$, and
$c_n\rho\le\tau\le\min\{\rho,1/3\}$; in particular $\tau>0$ and
$a_0\le1/2$. Then $a\in[a_0,1/2]$ and $a\le\sqrt{(1-a)/2}$ (which is
equivalent to $2a^2+a-1\le0$), so $b\in[a,\sqrt{(1-a)/2}]\subset[0,1]$.
Finally, $|P_i|^2=\tau+h^2=\rho$ for $i=1,2,3$.

\begin{proposition}[three-contact normal form]\label{prop:normal3}
Let $n\ge4$, and let $X\subset\R^n$ be compact with $\diam X=1$ and
circumradius greater than $1/2$. Then there are $u\in[0,1]^4$ and an
isometry $\varphi$ of $\R^n$ such that $\varphi(X)\subset H(u)$.
\end{proposition}

\begin{proof}
Translate the circumcentre to $0$ and let $\rho$ be the squared
circumradius; then $\rho\in(1/4,\rho_n]$ by Theorem~\ref{thm:jung}. Let
$P_1P_2P_3$ be the triangle of Lemma~\ref{lem:triangle}, labelled so that
its angles at $P_3$, $P_2$, $P_1$ are $\alpha\ge\beta\ge\gamma$. An
orthogonal map takes the circumcentre of the triangle to $he_3$ and its
plane to $he_3+\operatorname{span}(e_1,e_2)$; a rotation of
$\operatorname{span}(e_1,e_2)$ then moves $P_1$ to
$(\sqrt\tau,0,h,0,\dots,0)$, and, if necessary, the reflection
$x_2\mapsto-x_2$ makes the second coordinate of $P_2$ positive. Since the
triangle is acute, the arcs $P_1P_2$, $P_1P_3$, $P_2P_3$ of its circumcircle
not containing the third vertex have central angles $2\alpha$, $2\beta$,
$2\gamma$ by the inscribed angle theorem. Hence $P_2$ and $P_3$ are the
points of the circle at angles $2\alpha$ and $-2\beta$ from $P_1$, which
are the formulas of Definition~\ref{def:H} with $a=\cos\alpha$ and
$b=\cos\beta$.

It remains to find $u$. Choose $u_1$ with $\rho(u_1)=\rho$. Since
$\tau\le1/3$ and $\tau\le\rho$, we have
$h=\sqrt{\rho-\tau}\ge\sqrt{\max\{0,\rho-1/3\}}=\max\{t,0\}$, and
$\tau\ge c_n\rho$ gives $h\le\sqrt{(1-c_n)\rho}$; choose $u_2$ accordingly.
Since $60^\circ\le\alpha<90^\circ$, we have $0<a\le1/2$, and the side
$|P_1-P_2|=2\sqrt\tau\sin\alpha\le1$ gives $a^2\ge1-1/(4\tau)$, that is,
$a\ge a_0$; choose $u_3$. Finally, $\beta\le\alpha$ gives $b\ge a$, and
$\gamma\le\beta$ gives $\beta\ge90^\circ-\alpha/2$, that is,
$b\le\sin(\alpha/2)=\sqrt{(1-a)/2}$; choose $u_4$.

Since $P_1,P_2,P_3\in X$, $X\subset B(0,\sqrt\rho)$ and $\diam X=1$, we
have $X\subset B(0,\sqrt\rho)\cap\bigcap_iB(P_i,1)$. To place $X$ in $F_n$
without moving $0$ and the $P_i$, let
$W=\operatorname{span}(e_4,\dots,e_n)$ and let $M$ be the orthogonal
projection of $X$ onto $W$. Lemma~\ref{lem:frame} with $d=n-3$ gives an
orthonormal basis $v_1,\dots,v_{n-3}$ of $W$. The orthogonal map $\psi$
that sends $v_i$ to $e_{i+3}$ and is the identity on
$\operatorname{span}(e_1,e_2,e_3)$ fixes $0$ and the $P_i$, so
$\psi(X)\subset B(0,\sqrt\rho)\cap\bigcap_iB(P_i,1)$. For $x\in X$ and
$1\le i\le n-3$ we have
$\ip{\psi(x)}{e_{i+3}}=\ip x{v_i}\in[-h_M(-v_i),\,h_M(v_i)]$, and
therefore $\psi(X)\subset F_n$.
\end{proof}

If the circumradius of $X$ equals $1/2$, then $X$ lies in a ball of
diameter one, which Borsuk's theorem~\cite{Borsuk} splits into $n+1$ parts
of smaller diameter. This case needs no computation.

\subsection{From a continuous family to finitely many bodies}\label{sec:boxes}

We subdivide the parameter cube into dyadic boxes
$Q=\prod_{i=1}^4[k_i2^{-l_i},(k_i+1)2^{-l_i}]$. For each box the formulas
of Definition~\ref{def:H} are evaluated in interval
arithmetic~\cite{MooreKearfott}. The endpoints of all intervals lie on the
grid $2^{-44}\Z$ and are stored as integers, and every operation, including
the square root (computed by the integer square root), is rounded outward,
so the enclosures are rigorous. This interval class is taken from the
verification program of~\cite{ABPR}. Where a formula is monotone in a
parameter, the bounds are taken at the corresponding faces of the box to
avoid the dependency problem: $h$ increases in $u_1$ and decreases in
$u_2$; $a$ decreases in $u_3$ and increases in $a_0$; and, for fixed bounds
on $a$, the lower and upper bounds for $b$ are affine and decreasing in
$u_4$.

The output for a box $Q$ consists of rational numbers $\rho^-\le\rho^+$
with $\rho(u)\in[\rho^-,\rho^+]$ for $u\in Q$, and, for each $i$ and each
coordinate $j$, a rational interval $[\bar P_{ij}-e_{ij},\bar P_{ij}+e_{ij}]$
that contains the $j$-th coordinate of $P_i(u)$ for all $u\in Q$. Let
$\bar P_i=(\bar P_{ij})_j$, and let $E_i\ge0$ be a rational number with
$E_i^2\ge\rho^+\sum_je_{ij}^2$ (an integer square root rounded up). Put
\[
 \Omega_Q=B\bigl(0,\sqrt{\rho^+}\bigr)\cap\bigcap_{i=1}^3
 B\Bigl(\bar P_i,\sqrt{1-\rho^-+2E_i+|\bar P_i|^2}\Bigr)\cap F_n .
\]
All centres and squared radii here are rational.

\begin{lemma}\label{lem:omega}
If $u\in Q$, then $H(u)\subset\Omega_Q$.
\end{lemma}

\begin{proof}
Let $x\in H(u)$. Then $|x|^2\le\rho\le\rho^+$. From $|x-P_i|^2\le1$ and
$|P_i|^2=\rho$ we get $|x|^2-2\ip{P_i}x\le1-\rho$. Moreover,
$|\ip{P_i-\bar P_i}x|\le|P_i-\bar P_i|\,|x|\le\bigl(\sum_je_{ij}^2\bigr)^{1/2}
(\rho^+)^{1/2}\le E_i$. Hence $|x|^2-2\ip{\bar P_i}x\le1-\rho^-+2E_i$, which
is the inequality $|x-\bar P_i|^2\le1-\rho^-+2E_i+|\bar P_i|^2$.
\end{proof}

\begin{proposition}[reduction to finitely many bodies]\label{prop:reduction3}
Let $n\ge4$ and $k\ge n+1$, and let $\mathcal Q$ be a finite family of
dyadic boxes whose union is $[0,1]^4$. Suppose that for every
$Q\in\mathcal Q$ there are closed sets $C_1,\dots,C_k$ covering $\R^n$ with
$\diam(\Omega_Q\cap C_i)<1$ for all $i$. Then $b(n)\le k$.
\end{proposition}

\begin{proof}
Let $A\subset\R^n$ be bounded with $\diam A>0$; after a homothety,
$\diam A=1$. The closure $X=\overline A$ is compact with diameter one. If
its circumradius is $1/2$, then $X$ lies in a ball of diameter one, which
can be split into $n+1\le k$ parts of smaller diameter~\cite{Borsuk}.
Otherwise Proposition~\ref{prop:normal3} gives $u\in[0,1]^4$ and an isometry
$\varphi$ with $\varphi(X)\subset H(u)$. Choose $Q\in\mathcal Q$ containing
$u$; then $\varphi(X)\subset\Omega_Q$ by Lemma~\ref{lem:omega}. Take the
sets $C_1,\dots,C_k$ for $Q$ and put
\[
 Y_i=\{x\in A:\ \varphi(x)\in C_i,\ \varphi(x)\notin C_j\ \text{for}\ j<i\}.
\]
The sets $Y_i$ partition $A$, and since $\varphi$ is an isometry,
$\diam Y_i\le\diam(\Omega_Q\cap C_i)<1$.
\end{proof}

\subsection{Power diagrams and admissible partitions}\label{sec:power}

By Proposition~\ref{prop:reduction3}, to prove $b(n)\le k$ it suffices to
find, for each box of a dyadic subdivision of $[0,1]^4$, a covering of
$\R^n$ by $k$ closed sets that cut $\Omega_Q$ into parts of diameter less
than one. We take the cells of a power diagram and enclose each part in a
polytope with rational data.

Let $a_1,\dots,a_k\in\Q^n$ and $\beta_1,\dots,\beta_k\in\Q$, and put
\[
 C_i=\{x\in\R^n:\ \ip{a_i}x+\beta_i\ge\ip{a_j}x+\beta_j\ \text{for all } j\}.
\]
These are the cells of a power diagram~\cite{Aurenhammer}; for
$a_i=2z_i$ and $\beta_i=w_i-|z_i|^2$ one obtains the cells of the weighted
Voronoi diagram with centres $z_i$ and weights $w_i$. Each $C_i$ is a
closed convex polyhedron, and $\bigcup_iC_i=\R^n$ because every $x$ lies in
a cell where $\max_j(\ip{a_j}x+\beta_j)$ is attained. The covering property
therefore never has to be checked.

\begin{lemma}\label{lem:support}
Let $a\in\Q^n\setminus\{0\}$, $\beta\in\Q$, $c\in\Q^n$ and $s\in\Q$ with
$s>0$. Then $B(c,\sqrt s)\subset\{x:\ip ax\le\beta\}$ if and only if
$\beta-\ip ac\ge0$ and $(\beta-\ip ac)^2\ge s|a|^2$.
\end{lemma}

\begin{proof}
The maximum of $\ip ax$ over $B(c,\sqrt s)$ is $\ip ac+\sqrt s\,|a|$.
\end{proof}

Both conditions of Lemma~\ref{lem:support} can be checked in rational
arithmetic. We also use the elementary fact that the diameter of a nonempty bounded
polyhedron $K$ equals the largest distance between two of its vertices:
$K\times K$ is the convex hull of $\vertset K\times\vertset K$, and the
convex function $(x,y)\mapsto|x-y|$ attains its maximum over a polytope at
a vertex.

The body $\Omega_Q$ is the intersection of $F_n$ with four balls with
rational centres and rational squared radii, namely $B(0,\sqrt{\rho^+})$
and $B\bigl(\bar P_i,\sqrt{1-\rho^-+2E_i+|\bar P_i|^2}\bigr)$, $i=1,2,3$.

\begin{definition}\label{def:admissible}
Let $Q$ be a dyadic box. An \emph{admissible $k$-partition} for $Q$
consists of rational $a_i,\beta_i$ ($1\le i\le k$), which define the cells
$C_i$, and, for each $i$, of finitely many rational halfspaces
$\Pi_{i,1},\dots,\Pi_{i,s_i}$, each attached to one of the four balls
defining $\Omega_Q$, such that
\begin{enumerate}[label=\textup{(\roman*)}]
\item each $\Pi_{i,j}$ contains its ball;
\item the polyhedron $K_i=F_n\cap C_i\cap\Pi_{i,1}\cap\dots\cap\Pi_{i,s_i}$ is
bounded;
\item $|v-w|^2<1$ for all vertices $v,w$ of $K_i$.
\end{enumerate}
\end{definition}

\begin{figure}[t]
\centering
\includegraphics[width=0.5\textwidth]{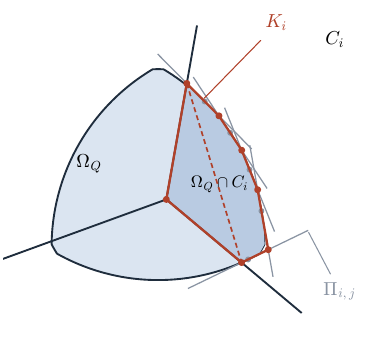}
\caption{A planar sketch of Definition~\ref{def:admissible}. The body
$\Omega_Q$ is an intersection of balls. For a cell $C_i$ of the power
diagram, the polygon $K_i$ is cut out by the walls of $C_i$ and by
halfplanes $\Pi_{i,j}$ tangent to the balls, so that
$K_i\supset\Omega_Q\cap C_i$. Its diameter is attained at two vertices
(dashed).}
\label{fig:admissible}
\end{figure}

Figure~\ref{fig:admissible} illustrates the definition in the plane.

\begin{proposition}\label{prop:admissible}
If $Q$ has an admissible $k$-partition, then
$\diam(\Omega_Q\cap C_i)\le\diam K_i<1$ for every $i$.
\end{proposition}

\begin{proof}
By (i), each $\Pi_{i,j}$ contains one of the balls defining $\Omega_Q$ and
hence $\Omega_Q$; since also $\Omega_Q\subset F_n$, we get
$\Omega_Q\cap C_i\subset K_i$. If $K_i$ is empty there is nothing to prove. Otherwise, by (ii), $K_i$ is
a polytope, so its
diameter is the largest distance between two vertices, which is less than
one by (iii).
\end{proof}

Only outer approximations enter the argument. The numerical estimates
used to find the partitions (Appendix~\ref{app:search}) do not enter their
verification.
Conditions (i)--(iii) are checked by the program described in
Appendix~\ref{app:verify}: (i) exactly by Lemma~\ref{lem:support}, (ii) by
an exact vertex enumeration (Proposition~\ref{prop:closure}), and (iii) in
exact arithmetic for every pair of vertices that a floating-point filter
with a proved error bound cannot settle.

\subsection{The partition trees and the end of the proof}\label{sec:trees}

\begin{figure}[ht]
\centering
\includegraphics[width=\textwidth]{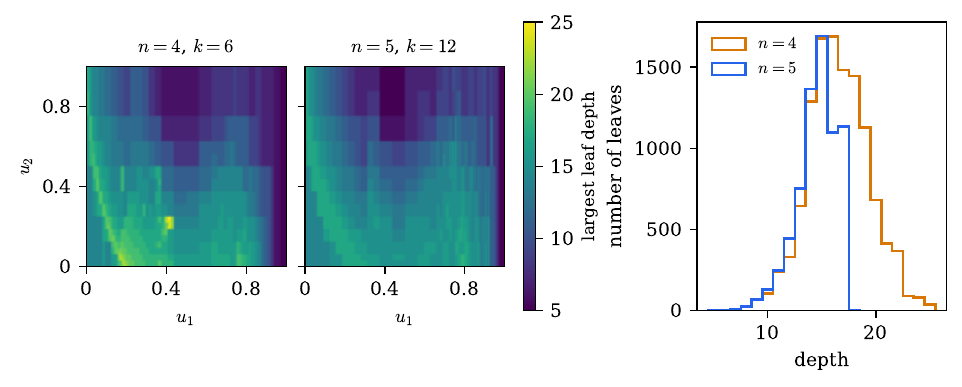}
\caption{The partition trees for $n=4$ and $n=5$. Left and centre: the
largest depth of a leaf over each point of the projection to the
$(u_1,u_2)$-plane. The parameter $u_1$ determines the circumradius;
$u_2=0$ corresponds to $\tau=c_n\rho$ and $u_2=1$ to
$\tau=\min\{\rho,1/3\}$. Right: the distribution of leaf depths.}
\label{fig:contact-tree}
\end{figure}

For each of the two dimensions the subdivision is stored as a binary tree
of dyadic boxes, and every leaf carries a power diagram with six or twelve
cells. Table~\ref{tab:contact} gives the sizes of the trees, and
Figure~\ref{fig:contact-tree} shows where the subdivision is fine.

\begin{table}[ht]
\centering
\caption{The partition trees for $n=4$ and $n=5$. The last column is an
upper bound for the diameters of all parts; it includes an outward
correction of the values printed in the verification logs
(Appendix~\ref{app:float}).}
\label{tab:contact}
\begin{tabular}{@{}rrrrrl@{}}
\toprule
$n$&Parts&Leaves&Depth&Diagrams&Diameter at most\\\midrule
4&6&11819&5--25&10615&$0.9999999994$\\
5&12&6976&5--18&5954&$0.99999995$\\\bottomrule
\end{tabular}
\end{table}

\begin{proof}[Proof of Theorem~\ref{thm:main45}]
The complete verification runs checked all $11819$ leaves for $n=4$ and
all $6976$ leaves for $n=5$, with no failures and covered volume one; the
largest printed diameters were $0.9999999993880884$ and
$0.9999999483201589$, respectively. Thus every leaf has an admissible
$6$-partition for $n=4$ and an admissible $12$-partition for $n=5$, and
Propositions~\ref{prop:admissible} and~\ref{prop:reduction3} give
$b(4)\le6$ and $b(5)\le12$.
\end{proof}

The numbers six and twelve are simply the numbers of cells in the diagrams
we found; we do not claim that they are optimal for the four-parameter
family. The diameter bounds in Table~\ref{tab:contact} are close to one
because a box is refined only until its polytopes pass the strict test.
They need not be close to the actual diameters of the parts of the curved
bodies $\Omega_Q$.

\section{Proof of Theorem~\ref{thm:main678}: coordinate vectors and a thirty-set covering}\label{sec:thm2}

By Proposition~\ref{prop:lens}, Theorem~\ref{thm:main678} follows from
coverings of $S^5$ and $S^4$ with small tolerance. The covering of $S^5$
comes from the normal regions of a weighted configuration and is treated
analytically (Section~\ref{sec:axis}). The covering of $S^4$ by thirty
sets (Section~\ref{sec:six}) is checked with a finite amount of rational
arithmetic; the matrices it uses are listed in Appendix~\ref{app:rays}.

\subsection{Coordinate vectors and the half-cube}\label{sec:axis}

For $m\ge6$ and $a>0$ consider the configuration
\begin{equation}\label{eq:axis-config}
 V^A_{m,a}=\{\pm ae_i:1\le i\le m\}\cup H_m
\end{equation}
of $2m+2^{m-1}$ vectors. It is invariant under coordinate permutations
and under changes of an even number of signs, and this group acts
transitively on $\{\pm ae_i\}$ and on $H_m$. It therefore suffices to study
the regions of $ae_1$ and of $h$.

\begin{lemma}\label{lem:axis-quadratic}
Let $k\ge5$ and $1\le q\le3/2$. If $1\ge b_1\ge\cdots\ge b_k\ge0$ and
$\sum_ib_i-2b_k\le q$, then $\sum_ib_i^2\le1+(q-1)^2$.
\end{lemma}

\begin{proof}
Put $d_i=b_i-b_{i+1}$ for $i<k$ and $d_k=b_k$. Then $b_i=\sum_{j\ge i}d_j$,
and the constraints read
\[
 d_i\ge0,\qquad \sum_i d_i\le1,\qquad
 \sum_{i<k}i\,d_i+(k-2)d_k\le q .
\]
The function $\sum_ib_i^2$ is a convex quadratic function of $d$, so by
Lemma~\ref{lem:vertices} (with weights $1,2,\dots,k-1,k-2$) it suffices to
evaluate it at the vertices of the form (a) and (b).

\emph{Vertices of type (a).} If only $d_i$ is nonzero, then
$\sum b_i^2=i\,d_i^2$. This gives $1$ for $i=1$, at most
$q^2/i\le q^2/2\le1+(q-1)^2$ for $2\le i<k$, and for $i=k$
\[
 \frac{kq^2}{(k-2)^2}\le\frac{5q^2}{9}\le1+(q-1)^2;
\]
the last inequality is equivalent to $2(2q-3)(q-3)\ge0$.

\emph{Vertices of type (b).} Both linear constraints are equalities, so
one of the two weights is less than $q<2$; hence one of the two nonzero
coordinates is $d_1$. The pairs $(d_1,d_j)$ with $j<k$ give the value
$1+(q-1)^2/(j-1)$, and the pair $(d_1,d_k)$ gives
$1+(k-1)(q-1)^2/(k-3)^2\le1+(q-1)^2$.
\end{proof}

\begin{proposition}\label{prop:axis}
Let $m\ge6$ and $2\le a\le5/2$, and put
\[
 \ell=\frac{m-3}{a-1},\qquad
 \eta_A=\frac{(a-2)^2}{2+(a-2)^2},\qquad
 \eta_H=\max\left\{0,\frac{2\ell-m+2}{\ell^2+m-1}\right\}.
\]
Then, with $\eta(m,a)=\max\{\eta_A,\eta_H\}$,
\[
 \min_{v\in V^A_{m,a}}\ \min_{x,y\in C_v(V^A_{m,a})}\ip xy=-\eta(m,a),
 \qquad\text{and hence}\qquad K(m,\eta(m,a))\le2m+2^{m-1}.
\]
\end{proposition}

\begin{proof}
\emph{The region of $ae_1$.} Normalize a point of the cone by $x_1=1$.
Comparison with $\pm ae_i$ gives $|x_i|\le1$. Let $b_1\ge\dots\ge b_{m-1}$
be the absolute values of the other coordinates. Choosing the signs of
$\varepsilon\in H_m$ to agree with those of $x$, except possibly at a
coordinate of smallest absolute value, we get
$\max_{\varepsilon\in H_m}\ip x\varepsilon\ge1+\sum b_i-2b_{m-1}$, so
comparison with $H_m$ gives $\sum b_i-2b_{m-1}\le a-1$.
Lemma~\ref{lem:axis-quadratic} with $k=m-1$ and $q=a-1$ yields
$|x|^2\le2+(a-2)^2$. Hence $\ip{x/|x|}{e_1}\ge(2+(a-2)^2)^{-1/2}$, and
Lemma~\ref{lem:caps} with $\gamma=1/(2+(a-2)^2)$ shows that all inner
products in the region are at least $2\gamma-1=-\eta_A$. Equality holds for
the normalizations of
\[
 (1,1,a-2,0,\ldots,0)\quad\text{and}\quad(1,-1,-(a-2),0,\ldots,0);
\]
both have absolute coordinate sum $a$ and largest absolute coordinate one,
so they belong to the cone.

\emph{The region of $h$.} The cone of this region is
\begin{equation}\label{eq:axis-sign}
 x_i+x_j\ge0\quad(i\ne j),\qquad a|x_i|\le S_x\quad(1\le i\le m).
\end{equation}
Here the pair inequalities express comparison with the vectors of $H_m$
that differ from $h$ in two coordinates; they imply the inequalities for
all even subsets, and hence comparison with all of $H_m$. Since $a\ge2$,
Lemma~\ref{lem:signcone} applies, and it remains to consider two points
whose negative coordinates are in different positions. We rescale these
coordinates to $-1$ and permute coordinates. Then $x$ and $y$ have the
form of Lemma~\ref{lem:negative-pair} with $k=m-2$, $c=a-1$ and
$\ell=(m-3)/(a-1)$: indeed, $ab\le S_x$ is equivalent to
$\sum p_i\ge(a-1)(b-\ell)$, and similarly for $y$. In the stated range
$1\le\ell<k$, and the lemma gives the bound $-\eta_H$.

If $\eta_H>0$, equality is attained by the normalizations of
\[
 (-1,\ell,1,\ldots,1)\quad\text{and}\quad(\ell,-1,1,\ldots,1),
\]
which satisfy~\eqref{eq:axis-sign} because $a\ell=\ell+m-3$. If
$\eta_H=0$, the bound for the coordinate regions is already the global
minimum.
\end{proof}

For $m=6$ and $a=7/3$ we obtain
\begin{equation}\label{eq:44}
 \eta_A=\frac1{19},\qquad \eta_H=\frac8{161}<\frac1{19},
 \qquad K\Bigl(6,\frac1{19}\Bigr)\le44.
\end{equation}
More generally, Proposition~\ref{prop:axis}, Lemma~\ref{lem:product} and
Proposition~\ref{prop:lens} give
\[
 b(n)\le1+2^{n-m-1}(2m+2^{m-1})
 \qquad\text{whenever }6\le m\le n-1\text{ and }\eta(m,a)<\muhat_n .
\]

\subsection{A covering of the four-sphere by thirty sets}\label{sec:six}

Let $\Sigma_k\subset\{-1,1\}^5$ denote the set of sign vectors with exactly
$k$ negative coordinates, and consider the configuration
\begin{equation}\label{eq:31-config}
 V=\{44e_i,\,-40e_i:1\le i\le5\}
 \cup20(\Sigma_0\cup\Sigma_1\cup\Sigma_2)\cup19\,\Sigma_4
\end{equation}
of $10+1+5+10+5=31$ vectors and its normal regions~\eqref{eq:normal-region}.

\begin{proposition}\label{prop:30}
The normal regions of $V$ can be modified to a covering of $S^4$ by $30$
closed sets in each of which all inner products are at least $-3/40$.
Consequently $K(5,3/40)\le30$.
\end{proposition}

The proof occupies the rest of this subsection.

\paragraph{Removing one region.}
Let $h_0=(1,\ldots,1)$ and $h_i=h_0-2e_i$, and write $C_0=C_{20h_0}(V)$ and
$C_i=C_{20h_i}(V)$ for $1\le i\le5$. Comparison with the vectors $20h_i$
shows that $C_0$ lies in the nonnegative orthant. On $C_i$, comparison
with $20h_0$ and with $20(h_i-2e_j)$ gives $x_i\le0$ and $x_j\ge0$ for
$j\ne i$, and comparison with $-40e_i$ gives $\sum_{j\ne i}x_j\ge-x_i$. We
replace $C_i$ by
\[
 \widetilde C_i=C_i\cup\{y\in C_0:y_i=\min_jy_j\},\qquad 1\le i\le5,
\]
and discard $C_0$. The new sets are closed and still cover $C_0$. All inner
products within $\widetilde C_i$ are nonnegative. For two points of $C_i$,
or two points of $C_0$, this follows from the sign conditions. If $x\in C_i$
and $y\in C_0$ with $y_i=\min_jy_j$, then $x_j\ge0$ and $y_j\ge y_i\ge0$
for $j\ne i$, whence
\[
 \ip xy\ge y_i\Bigl(x_i+\sum_{j\ne i}x_j\Bigr)\ge0.
\]

\paragraph{A bilinear criterion.}
Under coordinate permutations the remaining $25$ regions fall into four
classes, with representatives
\[
\begin{aligned}
 v_+&=(44,0,0,0,0),&v_-&=(-40,0,0,0,0),\\
 v_2&=(-20,-20,20,20,20),&v_4&=(-19,-19,-19,-19,19).
\end{aligned}
\]
The coordinates of each representative split into two blocks on which it
is constant, of sizes $(1,4)$, $(1,4)$, $(2,3)$ and $(4,1)$, respectively.
Every permutation of coordinates within the blocks preserves both $V$ and
the representative, and hence the region. We call a vector
\emph{increasingly ordered} if its coordinates increase within each block,
and denote by $P$ the permutation matrix that reverses the order within
each block.

\begin{lemma}\label{lem:metric}
Let $C\subset\R^5$ be a cone invariant under permutations within two fixed
blocks of coordinates. Let $S=I+E$ be a symmetric matrix that is invariant
under the same permutations and whose diagonal entries exceed the
off-diagonal entries within each block of size at least two. Suppose that
the increasingly ordered vectors of $C$ lie in a pointed polyhedral cone
$\Gamma=\{x:Ax\ge0\}$ with extreme rays generated by $r_1,\dots,r_N$. If
\[
 r_i^{\mathsf T}SP r_j\ge0\quad\text{for all }i,j
 \qquad\text{and}\qquad \|E\|_{\mathrm{op}}\le\mu,
\]
then $\ip xy\ge-\mu$ for all unit vectors $x,y\in C$.
\end{lemma}

\begin{proof}
By invariance, $S$ restricted to a block has the form $\alpha I+\beta J$
with $\alpha>0$, and its entries between the two blocks are constant. Hence
$x^{\mathsf T}Sy$ is the sum of a function of the block sums of $x$ and $y$
and of positive multiples of the inner products of the restrictions of $x$
and $y$ to the blocks. Applying a common permutation, we may assume that
$x$ is increasingly ordered. By the rearrangement
inequality~\cite[Theorem~368]{HLP}, $x^{\mathsf T}Sy$ does not increase
when $y$ is replaced by $Py^{\uparrow}$, where $y^{\uparrow}$ is the
increasingly ordered rearrangement of $y$. Both $x$ and $y^{\uparrow}$
belong to $C\cap\Gamma$, and expanding them in the generators $r_i$ gives
$x^{\mathsf T}SPy^{\uparrow}\ge0$. Thus $x^{\mathsf T}Sy\ge0$, and for unit
vectors $\ip xy=x^{\mathsf T}Sy-x^{\mathsf T}Ey\ge-\|E\|_{\mathrm{op}}$.
\end{proof}

\paragraph{Verification for the four classes.}
For each class we take a matrix $A$ whose rows are either positive
multiples of $v-w$ with $w\in V$, or ordering inequalities
$x_{j'}-x_j\ge0$ for consecutive coordinates $j,j'$ of a block. Thus the
increasingly ordered vectors of the cone over the region of $v$ lie in
$\Gamma=\{x:Ax\ge0\}$. The matrices $A_\pm,A_2,A_4$, the matrices
$E_\pm,E_2,E_4$, and the primitive integer generators of the extreme rays
of the four cones $\Gamma$ are listed in Appendix~\ref{app:rays}. The
hypotheses of Lemma~\ref{lem:metric} are checked as follows.

\begin{itemize}[leftmargin=*]
\item Each $A$ has rank five, so each $\Gamma$ is pointed. Every extreme
ray is the kernel of four linearly independent rows; it is obtained from
the signed $4\times4$ minors and kept if it satisfies $Ar\ge0$. There are
$\binom94+2\binom84+\binom74=301$ subsets of four rows to examine.
\item The norms of $E_+$ and $E_-$ are $13/200$ and $3/50$. The eigenvalues
of $E_2$ are $0$, $-3/50$, $-3/50$, $8/125$, $-3/125$, and those of $E_4$
are $-9/125$ (three times) and $(6\pm\sqrt{157})/250$. Since
$\sqrt{157}<51/4$, all four norms are less than $3/40$.
\item The differences between diagonal and off-diagonal entries of $S$ in
the nontrivial blocks are $187/200$, $47/50$, $1$, $47/50$ and $116/125$.
\item The minima of $r_i^{\mathsf T}SPr_j$ are given in
Table~\ref{tab:matrix}; there are $1007$ products in total.
\end{itemize}

\begin{table}[ht]
\centering
\caption{The four classes of regions in Section~\ref{sec:six}: sizes of the
data and the smallest value of the bilinear form.}
\label{tab:matrix}
\begin{tabular}{@{}crrrr@{}}
\toprule
Class&Rows of $A$&Rays&Pairs&$\min r_i^{\mathsf T}SP r_j$\\\midrule
$+$&9&17&289&$89/100$\\
$-$&8&15&225&$53/50$\\
$2$&8&18&324&$12/125$\\
$4$&7&13&169&$209/100$\\\bottomrule
\end{tabular}
\end{table}

All of this is a finite computation with integers and rationals, and the
norm bounds can be confirmed by Sylvester's criterion applied to
$\tfrac3{40}I\pm E$. A short program that performs these checks is
described in Appendix~\ref{app:rays}.

\begin{proof}[Proof of Proposition~\ref{prop:30}]
Lemma~\ref{lem:metric}, applied to the four classes, shows that each of
the $25$ regions other than $C_0,\dots,C_5$ has inner products at least
$-3/40$. Together with the five sets $\widetilde C_i$ this gives a covering
of $S^4$ by $30$ sets.
\end{proof}

\subsection{Proof of Theorem~\ref{thm:main678}}\label{sec:proof2}

\begin{proof}[Proof of Theorem~\ref{thm:main678}]
Since $D(6;3,40)=595>0$, Proposition~\ref{prop:30} and
Proposition~\ref{prop:lens} give $b(6)\le1+30=31$. By~\eqref{eq:44} and
$D(7;1,19)=680>0$, Proposition~\ref{prop:lens} gives $b(7)\le1+44=45$.
Finally, Lemma~\ref{lem:product} and the monotonicity of $K$ in $\mu$
give $K(7,1/19)\le K(6,1/19)\,K(1,1/19)\le44\cdot2=88$, and since
$D(8;1,19)=288>0$, we get $b(8)\le89$.
\end{proof}

\section{Proof of Theorem~\ref{thm:analytic}: the star configuration}\label{sec:thm3}

For $m\ge7$ consider the configuration
\begin{equation}\label{eq:star-config}
 V^S_m=\{2(\varepsilon e_i+\varepsilon' e_m):\ 1\le i<m,\
 \varepsilon,\varepsilon'\in\{-1,1\}\}\cup H_m
\end{equation}
of $N^S_m=4(m-1)+2^{m-1}$ vectors; thus $N^S_7=88$, $N^S_8=156$ and
$N^S_9=288$. Besides the half-cube it contains only those vectors
$2(\pm e_i\pm e_j)$ whose support contains the last coordinate. For $m=8$
these are $156$ of the $240$ roots of the lattice $E_8$, multiplied by two.
The configuration is invariant under permutations of the first $m-1$
coordinates and under changes of an even number of signs, and this group
acts transitively on each of the two parts of~\eqref{eq:star-config}. It
therefore suffices to study the regions of $2(e_1+e_m)$ and of $h$.

\subsection{The normal regions}\label{sec:star-regions}

\begin{lemma}\label{lem:bilinear}
Let $A,B,s,s'\ge0$, and let $p,q\in\R^k$ have nonnegative coordinates
with
\[
 p_i\le A,\qquad \sum_ip_i\le A+s,\qquad q_i\le B,\qquad \sum_iq_i\le B+s' .
\]
Then $\sum_ip_iq_i\le AB+ss'$.
\end{lemma}

\begin{proof}
If $s\ge A$, then $\sum p_iq_i\le A\sum q_i\le A(B+s')\le AB+ss'$. Let
$s<A$. Adding a zero coordinate if necessary, we may assume $k\ge2$, and by
the rearrangement inequality~\cite{HLP} we may assume
$p_1\ge\cdots\ge p_k$ and $q_1\ge\cdots\ge q_k$. Then
\[
 \sum_ip_iq_i\le p_1q_1+q_2\sum_{i\ge2}p_i\le p_1(q_1-q_2)+q_2(A+s)
 \le Aq_1+sq_2=s(q_1+q_2)+(A-s)q_1 ,
\]
and the right-hand side is at most $s(B+s')+(A-s)B=AB+ss'$.
\end{proof}

\begin{proposition}\label{prop:star}
Let $m\ge7$. All inner products within each normal region of $V^S_m$ are
at least
\[
 -\kappa_m,\qquad \kappa_m=\max\Bigl\{0,\frac{m-8}{(m-5)^2+m-1}\Bigr\}.
\]
Consequently
\[
 K(7,0)\le88,\qquad K(8,0)\le156,\qquad K(9,\tfrac1{24})\le288 .
\]
\end{proposition}

\begin{proof}
\emph{The region of $2(e_1+e_m)$.} We show that $\ip xy\ge0$ for any two
points $x,y$ of the cone over this region. Put $k=m-2\ge5$ and
$I=\{2,\dots,m-1\}$, and write
\[
 A=x_1,\quad s=x_m,\quad p_i=|x_i|\qquad(i\in I).
\]
Comparison with $2(-e_1+e_m)$ and $2(e_1-e_m)$ gives $A\ge0$ and $s\ge0$,
and comparison with $2(\pm e_i+e_m)$ gives $p_i\le A$. Call $x$
\emph{odd} if all its coordinates are nonzero and an odd number of them
are negative. If $x$ is not odd, some $\varepsilon\in H_m$ has
$\ip x\varepsilon=A+s+\sum p_i$. If $x$ is odd, we change the sign of
$\varepsilon$ at a coordinate of smallest absolute value; since
$p_i\le A$, this value is $\min\{s,\min_ip_i\}$. Comparison with $H_m$
therefore gives
\begin{equation}\label{eq:star-sum}
 \sum_{i\in I}p_i\le A+s+2\delta,\qquad
 \delta=\begin{cases}\min\{s,\min_ip_i\},&x\text{ odd},\\0,&\text{otherwise}.\end{cases}
\end{equation}
In both cases $0\le\delta\le s$ and $\delta\le p_i$ for all $i$. Define
$B,s',q_i$ and $\delta'$ in the same way for $y$, and let $D\subset I$ be
the set of indices with $x_iy_i<0$. Then
$\ip xy\ge AB+ss'-\sum_{i\in D}p_iq_i$, so it suffices to prove
\begin{equation}\label{eq:star-goal}
 \sum_{i\in D}p_iq_i\le AB+ss' .
\end{equation}
Put $r_i=p_i-\delta$ and $\sigma=s-(k-3)\delta$, and similarly
$r_i'=q_i-\delta'$ and $\sigma'=s'-(k-3)\delta'$. By~\eqref{eq:star-sum},
\begin{equation}\label{eq:star-r}
 0\le r_i\le A-\delta,\qquad \sum_{i\in I}r_i\le(A-\delta)+\sigma ,
\end{equation}
and the same holds for $r'$.

\emph{Case 1: $\delta=\delta'=0$.} Then~\eqref{eq:star-goal} is
Lemma~\ref{lem:bilinear}.

\emph{Case 2: $\delta>0=\delta'$.} We have
$\sum_Ip_iq_i=\delta\sum_Iq_i+\sum_Ir_iq_i$, and
$\delta\sum_Iq_i\le\delta(B+s')$. If $\sigma\ge0$, then
by~\eqref{eq:star-r} and Lemma~\ref{lem:bilinear} the second sum is at
most $(A-\delta)B+\sigma s'$, and altogether
$\sum_Ip_iq_i\le AB+ss'-(k-4)\delta s'\le AB+ss'$. If $\sigma<0$, the
second sum is at most $B\sum_Ir_i\le B(A-\delta)$, and
$\sum_Ip_iq_i\le AB+\delta s'\le AB+ss'$. The case $\delta'>0=\delta$ is
symmetric.

\emph{Case 3: $\delta,\delta'>0$.} Both points are odd. Since
$x_1,x_m,y_1,y_m$ are positive, each of $x$ and $y$ has an odd number of
negative coordinates in $I$, and therefore $d=|D|$ is even. If at least two
indices of $I$ lie outside $D$, then~\eqref{eq:star-sum} and
$p_i\ge\delta$ give $\sum_Dp_i\le A+s$, similarly $\sum_Dq_i\le B+s'$,
and~\eqref{eq:star-goal} follows from Lemma~\ref{lem:bilinear}. Let now
$d\ge k-1$; if $d=k$, then $k$ is even and so $k\ge6$. We have
\begin{equation}\label{eq:star-expand}
 \sum_{i\in D}p_iq_i=d\,\delta\delta'+\delta'\sum_Dr_i+\delta\sum_Dr_i'
 +\sum_Dr_ir_i' .
\end{equation}
Suppose first that $\sigma,\sigma'\ge0$. By~\eqref{eq:star-r} and
Lemma~\ref{lem:bilinear} the last sum is at most
$(A-\delta)(B-\delta')+\sigma\sigma'$, and~\eqref{eq:star-expand} is at
most
\[
 AB+(d-1)\delta\delta'+\delta'\sigma+\delta\sigma'+\sigma\sigma'.
\]
The difference between $AB+ss'$ and this expression is
\[
 (k-4)(\delta\sigma'+\delta'\sigma)
 +\bigl((k-5)(k-2)+k-d\bigr)\delta\delta'\ge0,
\]
since $k\ge5$, $d\le k$ and $\sigma,\sigma'\ge0$.
This proves~\eqref{eq:star-goal}.
Suppose now that $\sigma<0$. By~\eqref{eq:star-r} and
$r_i'\le B-\delta'$,
\[
 \delta'\sum_Dr_i+\sum_Dr_ir_i'\le B\sum_Dr_i\le B(A-\delta+\sigma)
 \le B(A-\delta)+\delta'\sigma ,
\]
because $B\ge\delta'$. With $\sum_Dr_i'\le B-\delta'+\sigma'$ we see
that~\eqref{eq:star-expand} is at most
\[
 AB+(d-1)\delta\delta'+\delta\sigma'+\delta'\sigma
 =AB+\delta s'+\delta's-(2k-5-d)\delta\delta'.
\]
Indeed the difference between $AB+ss'$ and the preceding expression is
\[
 (s-\delta)(s'-\delta')+(2k-6-d)\delta\delta'.
\]
The first term is nonnegative. Finally
$2k-6-d\ge0$: for $d=k-1$ this is $k\ge5$, and for $d=k$ it is $k\ge6$.
The case $\sigma'<0$ is symmetric.

\emph{The region of $h$.} Every point $x$ of the cone satisfies
\[
 x_i+x_j\ge0\quad(i\ne j),\qquad 2(|x_i|+|x_m|)\le S_x\quad(i<m);
\]
the first condition comes from comparison with the vectors of $H_m$ that
differ from $h$ in two coordinates, and the second from comparison with
$2(\pm e_i\pm e_m)$. In particular $2|x_i|\le S_x$ for all $i$, so by
Lemma~\ref{lem:signcone} it remains to consider two points $x,y$ whose
negative coordinates are in different positions $i\ne j$. Rescale so that
$x_i=y_j=-1$, and put $b=x_j$ and $d=y_i$; all coordinates of $x$ other
than $x_i$, and all coordinates of $y$ other than $y_j$, are at least $1$.
Write the remaining $k=m-2$ coordinates as $1+p_l$ and $1+q_l$. We claim
that
\begin{equation}\label{eq:star-h}
 \sum_lp_l\ge b-(m-5),\qquad \sum_lq_l\ge d-(m-5).
\end{equation}
If $j=m$, the condition $2(|x_i|+|x_m|)\le S_x$ reads
$2(1+b)\le-1+b+k+\sum p_l$, which is the first inequality. If $j<m$, we
use instead $2(|x_j|+|x_m|)\le S_x$ when $i<m$, that is,
$2(b+x_m)\le-1+b+k+\sum p_l$ with $x_m\ge1$, and
$2(|x_j|+|x_m|)=2(b+1)\le S_x$ when $i=m$. The second inequality
of~\eqref{eq:star-h} follows in the same way. Lemma~\ref{lem:negative-pair}
with $c=1$ and $\ell=m-5$ now gives the lower bound
$-\max\{0,(2\ell-k)/(\ell^2+k+1)\}=-\kappa_m$.

Finally $\kappa_7=\kappa_8=0$ and $\kappa_9=1/24$.
\end{proof}

The assumption $m\ge7$ cannot be dropped: for $m=6$ the points
$(1,1,1,1,-1,1)$ and $(1,-1,-1,-1,1,1)$ lie in the cone over the region of
$2(e_1+e_6)$, and their normalized inner product is $-1/3$.

\subsection{Proof of Theorem~\ref{thm:analytic}}\label{sec:proof3}

\begin{proof}[Proof of Theorem~\ref{thm:analytic}]
Let $n-1=7a+8b+c$. Applying Lemma~\ref{lem:product} repeatedly to
$K(7,0)\le88$, $K(8,0)\le156$ and $K(1,0)=2$, we get
$K(n-1,0)\le88^a156^b2^c$, and Proposition~\ref{prop:lens} with $\mu=0$
gives~\eqref{eq:star-bound}. The choice $a=0$, $b=\lfloor(n-1)/8\rfloor$
gives~\eqref{eq:eight}, because $156=2^8\cdot\tfrac{39}{64}$. Finally
$D(10;1,24)=495>0$, so $1/24<\muhat_{10}$ by~\eqref{eq:slack}, and
Proposition~\ref{prop:lens} with $K(9,\tfrac1{24})\le288$ gives
$b(10)\le289$.
\end{proof}

\section{Proof of Theorem~\ref{thm:ranges}: vectors supported on two coordinates}\label{sec:thm4}

For $m\ge9$ and a weight $\lambda>0$ consider the configuration
\begin{equation}\label{eq:root-config}
 V^D_{m,\lambda}=
 \{\lambda(\varepsilon_i e_i+\varepsilon_j e_j):i<j,
     \ \varepsilon_i,\varepsilon_j\in\{-1,1\}\}\cup H_m
\end{equation}
of $N_m=2m(m-1)+2^{m-1}$ vectors; it contains the configuration $V^S_m$
of Section~\ref{sec:thm3} with the weight $2$ replaced by $\lambda$. With
all vectors supported on two coordinates present, a weight slightly larger
than $2$ keeps the inner products within the regions close to zero also for
$m=11$ and $m=12$; we compute their minimum exactly
(Proposition~\ref{prop:roots}). The
configuration is invariant under coordinate permutations and under changes
of an even number of signs, and this group acts transitively on each of the
two parts of~\eqref{eq:root-config}; so it suffices to study the regions
of $\lambda(e_1+e_2)$ and of $h$.

\subsection{A quadratic inequality}

\begin{lemma}\label{lem:root-quadratic}
Let $k\ge7$, $2\le q\le12/5$, $a\ge0$ and $b_1\ge\cdots\ge b_k\ge0$. If
$a+b_1\le1$ and $\sum_ib_i-2b_k\le q$, then
\[
 2a^2+\sum_ib_i^2\le2+(q-2)^2.
\]
\end{lemma}

\begin{proof}
Put $d_i=b_i-b_{i+1}$ for $i<k$ and $d_k=b_k$. The constraints become
\[
 a,d_i\ge0,\qquad a+\sum_i d_i\le1,\qquad
 \sum_{i<k}i\,d_i+(k-2)d_k\le q,
\]
and the left-hand side of the claim is a convex quadratic function of
$(a,d_1,\dots,d_k)$. By Lemma~\ref{lem:vertices}, with the weights
$0,1,\dots,k-1,k-2$, it suffices to check the vertices of types (a)
and~(b).

\emph{Type (a).} The variables $a,d_1,d_2$ give the values $2,1,2$; the
variable $d_i$ with $3\le i<k$ gives $q^2/i\le48/25$; and $d_k$ gives
$kq^2/(k-2)^2\le1008/625$.

\emph{Type (b).} Both linear constraints are equalities, so $q$ lies
strictly between the two weights. The smaller weight is therefore $0$, $1$
or $2$, and the larger one is at least $3$.

For the pair $(a,d_i)$ with $3\le i<k$ we have $d_i=q/i$, and the value is
\[
 2\Bigl(1-\frac qi\Bigr)^2+i\Bigl(\frac qi\Bigr)^2
 =2+\frac q{i^2}\bigl((i+2)q-4i\bigr)\le2,
\]
because $4i/(i+2)\ge12/5\ge q$. For the pair $(a,d_k)$ the value is
$2(1-q/(k-2))^2+k(q/(k-2))^2$. It is at most two for $k\ge8$, because then
$q\le4(k-2)/(k+2)$. For $k=7$ it is smaller than $2+(q-2)^2$ by
$16(q-5/2)^2/25\ge0$.

For a pair $(d_i,d_j)$ with $i\in\{1,2\}$ and $j\ge3$, if $j<k$, the value
is
\[
 i+\frac{(q-i)^2}{j-i}\le
 \begin{cases}1+(q-1)^2/2,&i=1,\\2+(q-2)^2,&i=2,
 \end{cases}
\]
and $2+(q-2)^2$ exceeds the first bound by $(q-3)^2/2$. If $j=k$, the value
is $i+(k-i)(q-i)^2/(k-2-i)^2$; for $i=1$ it is at most
$1+3(q-1)^2/8\le347/200<2$, and for $i=2$ it is at most
$2+5(q-2)^2/9$.
\end{proof}

\subsection{The normal regions}

\begin{proposition}\label{prop:roots}
Let $m\ge9$ and $2\le\lambda\le11/5$, and put
\[
 \ell=\frac{m-3-\lambda}{\lambda-1},\qquad
 \mu_D=\frac{(\lambda-2)^2}{1+(\lambda-2)^2},\qquad
 \mu_H=\max\left\{0,\frac{2\ell-m+2}{\ell^2+m-1}\right\}.
\]
Then, with $\kappa(m,\lambda)=\max\{\mu_D,\mu_H\}$,
\[
 \min_{v\in V^D_{m,\lambda}}\ \min_{x,y\in C_v(V^D_{m,\lambda})}\ip xy
 =-\kappa(m,\lambda),\qquad\text{and hence}\qquad K(m,\kappa(m,\lambda))\le N_m.
\]
\end{proposition}

\begin{proof}
\emph{The region of $\lambda(e_1+e_2)$.} Comparison with the other vectors
supported on two coordinates gives $x_1+x_2\ge|x_i|+|x_j|$ for all
$i\ne j$; hence $x_1,x_2\ge0$ and $|x_j|\le\min\{x_1,x_2\}$ for $j\ge3$.
Scale so that $x_1+x_2=2$ and write $\{x_1,x_2\}=\{1+a,1-a\}$ with
$0\le a\le1$. Let $b_1\ge\cdots\ge b_k$, $k=m-2$, be the absolute values of
the remaining coordinates. Then $a+b_1\le1$. Choosing the signs of
$\varepsilon\in H_m$ to agree with those of $x$, except possibly at the
coordinate with absolute value $b_k$, we get
$\max_{\varepsilon\in H_m}\ip x\varepsilon\ge2+\sum b_i-2b_k$, so
comparison with $H_m$ gives $\sum b_i-2b_k\le2\lambda-2$. Lemma~\ref{lem:root-quadratic} with
$q=2\lambda-2$ yields $|x|^2=2+2a^2+\sum b_i^2\le4+4(\lambda-2)^2$. Since
$\ip x{(e_1+e_2)/\sqrt2}=\sqrt2$, Lemma~\ref{lem:caps} with
$\gamma=1/(2+2(\lambda-2)^2)$ gives the lower bound $2\gamma-1=-\mu_D$.

\emph{The region of $h$.} The cone of this region is
\begin{equation}\label{eq:root-sign}
 x_i+x_j\ge0,\qquad \lambda(|x_i|+|x_j|)\le S_x\qquad(i\ne j).
\end{equation}
Since $\lambda\ge2$, Lemma~\ref{lem:signcone} applies, and it remains to
consider two points whose negative coordinates are in different positions.
After rescaling and permuting coordinates they have the form of
Lemma~\ref{lem:negative-pair} with $k=m-2$, $c=\lambda-1$ and
$\ell=(k-1-\lambda)/(\lambda-1)$; the required bound for $\sum p_i$ follows
from $\lambda(b+1)\le S_x$. In the stated range of parameters
$1\le\ell<k$. This gives the lower bound $-\mu_H$.

\emph{Equality.} In the first type of region take
\[
 (1,1,1,1,2\lambda-4,0,\ldots,0)\quad\text{and}\quad
 (1,1,-1,-1,-(2\lambda-4),0,\ldots,0).
\]
For each of them the largest sum of two absolute coordinates is $2$ and the
absolute coordinate sum is $2\lambda$, so both lie in the cone; after
normalization their inner product is $-\mu_D$. In the region of $h$ take
$(-1,\ell,1,\ldots,1)$ and $(\ell,-1,1,\ldots,1)$. They
satisfy~\eqref{eq:root-sign} because $\ell+m-3=\lambda(\ell+1)$, and their
normalized inner product is $(m-2-2\ell)/(\ell^2+m-1)$. This equals
$-\mu_H$ if $\mu_H>0$; otherwise the first pair already attains the global
minimum.
\end{proof}

\subsection{Proof of Theorem~\ref{thm:ranges}}\label{sec:proof4}

\begin{proof}[Proof of Theorem~\ref{thm:ranges}]
We apply Proposition~\ref{prop:roots} with the following parameters:
\begin{center}
\begin{tabular}{@{}rrcrr@{}}
\toprule
$m$&$\lambda$&$(\mu_D,\mu_H)$&$\mu=\kappa(m,\lambda)$&$N_m$\\\midrule
11&$13/6$&$(1/37,\ 1/35)$&$1/35$&1244\\
12&$761/348$&$(4225/125329,\ 7021/208275)$&$4225/125329$&2312\\\bottomrule
\end{tabular}
\end{center}
This gives $K(11,\tfrac1{35})\le1244$ and
$K(12,\tfrac{4225}{125329})\le2312$. By Lemma~\ref{lem:product} with
one-dimensional factors, $K(n-1,\tfrac1{35})\le2^{n-12}\cdot1244$ for
$n\ge12$, and Proposition~\ref{prop:lens} gives $b(n)\le1+2^{n-12}\cdot1244$
whenever $1/35<\muhat_n$. Since $\muhat_n$ decreases in $n$ and
$D(14;1,35)=1800>0$, this holds for $n\le14$ by~\eqref{eq:slack}, which
proves $b(12)\le1245$ and $b(14)\le4977$. Similarly
$D(13;4225,125329)=1296974>0$ gives $b(13)\le2313$.
\end{proof}

\subsection{The values in Table~\ref{tab:bounds}}\label{sec:table}

Table~\ref{tab:blocks} lists all blocks constructed in
Sections~\ref{sec:thm2}--\ref{sec:thm4}. The entries of
Table~\ref{tab:bounds} for $6\le n\le19$ are the values of the
bound~\eqref{eq:block-formula} for these blocks. In dimensions $12$, $13$
and $14$ the bound comes from the blocks of dimensions $11$ and $12$, and
in dimension $10$ from the nine-dimensional block. All other entries with
$n\ge9$ are instances of~\eqref{eq:star-bound}: the exponents $(a,b,c)$ are
$(0,1,0)$, $(0,1,2)$, $(2,0,0)$, $(1,1,0)$, $(0,2,0)$, $(0,2,1)$ and
$(0,2,2)$ for $n=9,11,15,16,17,18,19$.

\begin{table}[tbp]
\centering
\caption{All blocks $(m,N,\mu)$ with $K(m,\mu)\le N$ used in this paper.
The weight is $a$ for the block from Proposition~\ref{prop:axis} and
$\lambda$ for the blocks from Proposition~\ref{prop:roots}.}
\label{tab:blocks}
\begin{tabular}{@{}rrrll@{}}
\toprule
$m$&$N$&$\mu$&Weight&Source\\\midrule
1&2&$0$&---&two antipodal points\\
5&30&$3/40$&---&Proposition~\ref{prop:30}\\
6&44&$1/19$&$7/3$&Proposition~\ref{prop:axis}\\
7&88&$0$&---&Proposition~\ref{prop:star}\\
8&156&$0$&---&Proposition~\ref{prop:star}\\
9&288&$1/24$&---&Proposition~\ref{prop:star}\\
11&1244&$1/35$&$13/6$&Proposition~\ref{prop:roots}\\
12&2312&$4225/125329$&$761/348$&Proposition~\ref{prop:roots}\\\bottomrule
\end{tabular}
\end{table}

\section{Concluding remarks}\label{sec:remarks}

\begin{remark}[uniform bounds]\label{rem:dnk}
All diameter bounds in our proofs are uniform. In Proposition~\ref{prop:lens}
every part has diameter at most $\Delta(n,\mu)<1$ by
Remark~\ref{rem:lens}. In Proposition~\ref{prop:reduction3} with the partitions
of Section~\ref{sec:trees}, every part lies in one of finitely many
polytopes $K_i$, whose diameters are bounded as in
Table~\ref{tab:contact}. Finally, a ball of diameter one is split into
$n+1$ fixed parts of diameter less than one. Consequently $d_{n,k}<1$ for
every pair $(n,k)$ for which Theorems~\ref{thm:main45}--\ref{thm:ranges}
give $b(n)\le k$; for example, $d_{4,6}<1$ and $d_{8,89}<1$.
\end{remark}

\paragraph{Large dimensions.} The factor $(39/64)^{\lfloor(n-1)/8\rfloor}$
in~\eqref{eq:eight} decays exponentially, but the base $156^{1/8}\approx1.88$
is much larger than $\sqrt{3/2}$. Thus~\eqref{eq:eight} is an improvement of Lassak's estimate,
but not of the best known asymptotic bound.

Better spherical
coverings alone cannot close this gap. Lassak's lens is a universal cover,
and by a theorem of Arman, Bondarenko, Prymak and
Radchenko~\cite{ABPRUniversal}, every partition of a universal cover in
$\R^n$ into sets of diameter at most one has at least $(\sqrt2-o(1))^n$
parts.

\paragraph{The spherical problem.} Proposition~\ref{prop:lens} reduces the
problem to estimating $K(m,\mu)$ for small positive $\mu$, and the threshold
$\muhat_n$ is sharp for this reduction (Remark~\ref{rem:lens}). The
configurations of Sections~\ref{sec:thm2}--\ref{sec:thm4} are only
examples. Already $K(8,0)\le156<2^8$ shows that orthants are not optimal,
and further improvements in $K(m,\mu)$ may yield sharper bounds in the
dimensions where the corresponding blocks are admissible, through
Corollary~\ref{cor:blocks}.

\paragraph{Small dimensions.} The reduction of Section~\ref{sec:thm1}
works in every dimension $n\ge4$, but the verification becomes
considerably heavier as $n$ grows: the complete run took about $75$
minutes for $n=4$ and about $11.5$ hours for $n=5$
(Appendix~\ref{app:runs}). Whether $b(4)=5$ remains open.

\section*{Acknowledgements}
The author thank Vsevolod Voronov for valuable discussions on the subject
of this paper.

\section*{Data availability}
The partition trees, the data for the covering of $S^4$, the verification
programs and the logs of the complete verification runs are available in
the repository~\cite{Repo}.

\section*{AI usage}
OpenAI's model GPT-6 Astra and Anthropic's model Claude Fable~5.2 were
used in obtaining the results of this paper, and Anthropic's model Claude
Opus~5.5 was used to edit the text of the manuscript.

\clearpage
\appendix
\section{Data for Proposition~\ref{prop:30}}\label{app:rays}

For each of the four classes of Section~\ref{sec:six}, the matrix $A$
defines the cone $\Gamma=\{x:Ax\ge0\}$, and $S=I+E$ is the matrix of
Lemma~\ref{lem:metric}. In classes $+$ and $-$ the permutation $P$ reverses
coordinates $2,3,4,5$; in class $2$ it swaps coordinates $1,2$ and $3,5$;
in class $4$ it reverses coordinates $1,2,3,4$.

\subsection*{The matrices $A$ and $E$}
{\small
\[A_{+}=\begin{pmatrix}0&-1&1&0&0\\0&0&-1&1&0\\0&0&0&-1&1\\1&0&0&0&-1\\6&-5&-5&-5&-5\\6&5&-5&-5&-5\\6&5&5&-5&-5\\11&10&0&0&0\\25&19&19&19&19\end{pmatrix},\qquad E_{+}=\frac1{200}\begin{pmatrix}13&0&0&0&0\\0&-13&0&0&0\\0&0&-13&0&0\\0&0&0&-13&0\\0&0&0&0&-13\end{pmatrix}.\]
\[A_{-}=\begin{pmatrix}-21&19&19&19&-19\\-10&0&0&0&-11\\-1&-1&-1&-1&-1\\-1&1&-1&-1&-1\\-1&1&0&0&0\\0&-1&1&0&0\\0&0&-1&1&0\\0&0&0&-1&1\end{pmatrix},\qquad E_{-}=\frac1{50}\begin{pmatrix}3&0&0&0&0\\0&-3&0&0&0\\0&0&-3&0&0\\0&0&0&-3&0\\0&0&0&0&-3\end{pmatrix}.\]
\[A_{2}=\begin{pmatrix}-5&-5&5&5&-6\\-1&-1&39&39&1\\-1&1&0&0&0\\0&-1&0&0&0\\0&-1&1&0&0\\0&0&-1&1&0\\0&0&0&-1&1\\1&-1&1&1&1\end{pmatrix},\qquad E_{2}=\frac1{250}\begin{pmatrix}5&5&-4&-4&-4\\5&5&-4&-4&-4\\-4&-4&-10&5&5\\-4&-4&5&-10&5\\-4&-4&5&5&-10\end{pmatrix}.\]
\[A_{4}=\begin{pmatrix}-19&-19&-19&-19&-25\\-1&1&0&0&0\\0&-1&1&0&0\\0&0&-1&1&0\\0&0&0&-1&1\\1&1&-39&-39&-1\\21&-19&-19&-19&19\end{pmatrix},\qquad E_{4}=\frac1{500}\begin{pmatrix}-21&15&15&15&-11\\15&-21&15&15&-11\\15&15&-21&15&-11\\15&15&15&-21&-11\\-11&-11&-11&-11&0\end{pmatrix}.\]
}

\subsection*{Extreme rays}
The rows of the following matrices are the primitive integer generators of
the extreme rays of the four cones $\Gamma$.
\begingroup\small\renewcommand{\arraystretch}{0.9}\setlength{\arraycolsep}{4pt}
\[R_{+}=\begin{pmatrix}5&-1&-1&-1&5\\5&-1&0&0&5\\5&0&0&1&5\\5&0&0&3&3\\5&0&2&2&2\\10&-11&0&0&1\\10&0&1&1&10\\10&3&3&3&3\\15&1&1&1&15\\20&-22&0&1&1\\30&-33&1&1&1\\76&-25&-25&-25&-25\\190&-209&-30&-30&19\\380&-418&-60&-11&-11\\380&-239&-239&-239&217\\380&-239&-239&-11&-11\\570&-627&-41&-41&-41\end{pmatrix},\qquad
R_{-}=\begin{pmatrix}-627&-41&-41&-41&570\\-418&-60&-11&-11&380\\-209&-30&-30&19&190\\-38&-38&-1&-1&2\\-38&-21&-21&-21&-21\\-38&-21&-21&19&19\\-33&1&1&1&30\\-22&0&1&1&20\\-19&-19&-2&-2&-2\\-19&-19&-2&1&1\\-11&0&0&1&10\\-4&1&1&1&1\\-3&0&1&1&1\\-2&0&0&1&1\\-1&-1&0&0&0\end{pmatrix}.\]
\[R_{2}=\begin{pmatrix}-418&-60&-11&-11&380\\-239&-239&-239&217&380\\-239&-239&-11&-11&380\\-209&-30&-30&19&190\\-40&-40&-40&37&37\\-38&-38&-1&-1&2\\-38&-21&-21&19&19\\-22&0&1&1&20\\-19&-19&-2&1&1\\-11&0&0&1&10\\-3&0&1&1&1\\-2&0&0&1&1\\-1&-1&-1&5&5\\-1&-1&0&0&0\\-1&0&0&5&5\\0&0&1&1&1\\0&0&1&5&5\\0&0&3&3&5\end{pmatrix},\qquad
R_{4}=\begin{pmatrix}-627&-41&-41&-41&570\\-418&-60&-11&-11&380\\-239&-239&-239&217&380\\-239&-239&-11&-11&380\\-209&-30&-30&19&190\\-40&-40&-40&37&37\\-38&-38&-1&-1&2\\-38&-21&-21&-21&-21\\-38&-21&-21&19&19\\-25&-25&-25&-25&76\\-19&-19&-2&-2&-2\\-19&-19&-2&1&1\\-1&-1&-1&-1&-1\end{pmatrix}.\]
\endgroup

\subsection*{Checking the data}
The program \path{checks/six.py} in the repository~\cite{Repo} reads these
matrices and ray lists from \path{data/six.json} and checks in integer and
rational arithmetic that every row of each $A$ has the origin claimed in
Section~\ref{sec:six}, that each $A$ has rank five, that the lists above
are exactly the extreme rays (by examining all $301$ subsets of four
rows), that each $S$ is invariant under the block permutations, that all
$1007$ products $r_i^{\mathsf T}SPr_j$ are nonnegative, and that
$\|E\|_{\mathrm{op}}<3/40$ (by Sylvester's criterion for
$\tfrac3{40}I\pm E$). The program \path{checks/analytic.py} rechecks the
arithmetic of Sections~\ref{sec:thm3} and~\ref{sec:thm4}: the block
parameters and their admissibility, the algebraic identities in the proof
of Proposition~\ref{prop:star}, and the entries of
Table~\ref{tab:bounds}. These checks supplement the proofs and do not
replace any part of them. Both programs use only the Python standard
library and are run by \texttt{python3 verify.py short}
(Appendix~\ref{app:repro}).

\section{Verification in dimensions four and five}\label{app:verify}

This appendix explains how conditions (i)--(iii) of
Definition~\ref{def:admissible} are checked for all leaves of the two
partition trees, and which parts of the computation are exact. The
partition trees, the verification code and the logs of the complete runs
are contained in the repository~\cite{Repo}, available at \repo; all
file names below refer to this repository.

\subsection{Overview}\label{app:overview}

For each dimension the partition tree is stored in preorder. Every internal
node splits one coordinate of its box at the midpoint, and every leaf
carries an index into a list of power diagrams, stored as integers over the
common denominator $10^9$. Apart from the floating-point steps named below,
the verification program works in exact arithmetic (Python's
arbitrary-precision integers and the class \texttt{Fraction}).

The program first reads the tree and rejects any unresolved leaf; by
induction on the tree, the leaves form a family of dyadic boxes with union
$[0,1]^4$, and as an additional check their volumes are added up exactly,
which gives~$1$. For every leaf $Q$ it then recomputes $\Omega_Q$ from the
endpoints of $Q$ by the interval arithmetic of Section~\ref{sec:boxes}, so
no geometric data are read from the stored tree. The tangent halfspaces
are regenerated by the procedure of Appendix~\ref{app:search} and checked
exactly, so how they were produced is irrelevant. Finally, each cell is
checked as described in Appendix~\ref{app:cellcheck}, using the exact
vertex enumeration of Appendix~\ref{app:enum} and the floating-point
filters of Appendix~\ref{app:float}.

\subsection{Exact vertex enumeration}\label{app:enum}

There are $70\,914$ polytopes $K_i$ for $n=4$ and $83\,712$ for $n=5$.
Vertex enumeration by the double description
method~\cite{FukudaProdon} in rational arithmetic is reliable but slow for
this many polytopes. We therefore let a fast floating-point convex hull
program (Qhull~\cite{Qhull}) propose the vertices and prove the
completeness of the proposed list by an exact adjacency argument, related
to the pivoting methods of~\cite{AvisFukuda}.

Let $K=\{x\in\R^n:\ Gx\le g\}$ with $G\in\Q^{m\times n}$ and $g\in\Q^m$,
and denote the rows of $G$ by $G_j$. For an index set $I$ of cardinality
$n$ such that the submatrix $G_I$ is nonsingular, put $x_I=G_I^{-1}g_I$. We
call $I$ \emph{confirmed} if $G_jx_I<g_j$ for all $j\notin I$. In this case
$x_I\in K$, exactly $n$ constraints are active at $x_I$, and $x_I$ is a
vertex of $K$.

\begin{proposition}[adjacency closure]\label{prop:closure}
Let $\mathcal S$ be a nonempty family of confirmed index sets such that for
every $I\in\mathcal S$ and every $i\in I$ there is $J\in\mathcal S$,
$J\ne I$, with $I\setminus\{i\}\subset J$. Then $K$ is bounded and
$\vertset K=\{x_I:\ I\in\mathcal S\}$.
\end{proposition}

\begin{proof}
Let $I\in\mathcal S$ and $x=x_I$. For $i\in I$ define $d_i$ by
$G_Id_i=-e_i$ (coordinates indexed by $I$). Moving from $x$ along $d_i$,
constraint $i$ becomes inactive and the other constraints of $I$ stay
active. The cone $T=\{d:\ G_Id\le0\}$ is generated by the $d_i$, since
every $d\in T$ equals $\sum_i(-(G_Id)_i)\,d_i$. All constraints outside
$I$ are strict at $x$, so $K\subset x+T$ and $x+td_i\in K$ for small
$t>0$. The line $\{y:G_{I\setminus i}y=g_{I\setminus i}\}$ is
$\{x+td_i:\ t\in\R\}$, and constraint $i$ is violated for $t<0$. Hence
$F_i=K\cap\{y:G_{I\setminus i}y=g_{I\setminus i}\}$ is a segment or a ray
with endpoint~$x$.

\emph{Step 1: each $F_i$ is a segment whose other endpoint is $x_J$ for
some $J\in\mathcal S$.} Take $J$ as in the hypothesis. Then
$G_{I\setminus i}x_J=g_{I\setminus i}$ and $x_J\in K$, so $x_J\in F_i$.
Moreover $x_J\ne x$, since otherwise more than $n$ constraints would be
active at $x$. Thus $x_J=x+t^*d_i$ with $t^*>0$. If $x+td_i\in K$ for some
$t>t^*$, the vertex $x_J$ would be an interior point of a segment in $K$,
which is impossible. Hence $F_i=[x,x_J]$.

\emph{Step 2: a monotone walk.} Let $c\in\R^n$ and $I_0\in\mathcal S$. If
at $x_{I_s}$ there is an $i$ with $\ip c{d_i}>0$, move along the edge $F_i$
to the vertex $x_{I_{s+1}}$ given by Step~1; then
$\ip c{x_{I_{s+1}}}>\ip c{x_{I_s}}$. Since $\mathcal S$ is finite and the
values increase strictly, the walk stops at some $x_I$ with
$\ip c{d_i}\le0$ for all $i$. Then $\ip cd\le0$ on $T$, and
$K\subset x_I+T$ shows that $x_I$ maximizes $\ip c\cdot$ on $K$.

\emph{Boundedness.} If $K$ were unbounded, it would contain a ray
$x_{I_0}+td$, $t\ge0$, with $d\ne0$. For $c=d$ the function $\ip c\cdot$ is
unbounded above on $K$, while by Step~2 it attains its maximum.

\emph{Completeness.} Let $w$ be a vertex of $K$ and $I(w)$ the set of
constraints active at $w$. The matrix $G_{I(w)}$ has rank $n$, for
otherwise $w\pm\eps d\in K$ for some $d\ne0$ and small $\eps>0$. Put
$c=\sum_{j\in I(w)}G_j$. For $y\in K$ we have
$\ip cy\le\sum_{j\in I(w)}g_j=\ip cw$, with equality only if
$G_{I(w)}y=g_{I(w)}$, i.e.\ $y=w$. So $w$ is the unique maximizer of
$\ip c\cdot$ on $K$, and by Step~2 it equals $x_I$ for some
$I\in\mathcal S$. Conversely, every $x_I$ with confirmed $I$ is a vertex.
\end{proof}

In practice $\mathcal S$ is the family of active sets returned by Qhull.
The hypothesis of Proposition~\ref{prop:closure} is checked combinatorially:
each $(n-1)$-element subset $I\setminus\{i\}$ must be contained in exactly
two sets of the family. If anything fails (a degenerate vertex, a missing
neighbour, a Qhull error, an empty or lower-dimensional $K$), the cell is
handled by \texttt{cddlib}~\cite{cdd} in GMP rational arithmetic.

\paragraph{Integer arithmetic.} Each row $(G_j,g_j)$ is multiplied by a
common denominator, which does not change $K$. The solution $x_I=X/D$ with
$X\in\Z^n$ and $D>0$ is computed by Bareiss' fraction-free
elimination~\cite{Bareiss}. To make the result independent of the
implementation of this step, two facts are checked separately: $G_IX=Dg_I$
in integers, and $\det G_I\not\equiv0\pmod{2^{61}-1}$. The latter implies
$\det G_I\ne0$, so $x_I$ is the unique solution of the system.

\subsection{Floating-point filters}\label{app:float}

To avoid tens of millions of comparisons in rational arithmetic, two
filters in IEEE~754 double precision are used. Let $\uround=2^{-53}$ be the
unit roundoff and $\gamma_k=k\uround/(1-k\uround)$. We use the standard
bound~\cite[\S3.1]{Higham}: if $\hat x,\hat y\in\R^k$ have floating-point
entries, then, in the absence of underflow and overflow, the computed inner
product differs from $\ip{\hat x}{\hat y}$ by at most
$\gamma_k\sum_t|\hat x_t||\hat y_t|$, for any order of summation. Rounding
a rational number $z$ to the nearest normal floating-point number has
relative error at most $\uround$. We assume round-to-nearest with gradual
underflow and finite intermediate results; subnormal numbers are covered by
the bound $|\operatorname{fl}(z)-z|\le\uround|z|+\eta$ with
$\eta=2^{-1075}$.

\begin{lemma}[strict inactivity filter]\label{lem:filter1}
Let $s=g_j-G_jx$ with rational $g_j,G_j,x$, where $n\le8$, and let $\hat s$
be the value computed from the rounded data $\hat g_j,\hat G_j,\hat x$. Put
$\vartheta=10^{-12}\bigl(|\hat g_j|+\sum_t|\hat G_{jt}||\hat x_t|+1\bigr)$.
If $\hat s>\vartheta$, then $s>0$.
\end{lemma}

\begin{proof}
Write $\Sigma=|g_j|+\sum_t|G_{jt}||x_t|$. Each of $\hat g_j,\hat G_{jt},
\hat x_t$ differs from the exact value by a factor of at most
$1+\uround$, so $|\hat g_j-g_j|\le\uround|g_j|$ and
$|\hat G_{jt}\hat x_t-G_{jt}x_t|\le\gamma_2|G_{jt}x_t|$. Evaluating
$\hat g_j-\sum_t\hat G_{jt}\hat x_t$, an inner product of length $n+1$,
adds at most $\gamma_{n+1}(1+\gamma_2)\Sigma$. Since
$\gamma_a+\gamma_b+\gamma_a\gamma_b\le\gamma_{a+b}$~\cite[Lemma~3.3]{Higham},
we get $|\hat s-s|\le\gamma_{n+3}\Sigma<1.3\cdot10^{-15}\,\Sigma$ for
$n\le8$. Similarly, the computed value of $\vartheta$ is at least
$10^{-12}(1-\gamma_{n+6})(\Sigma+1)>0.9\cdot10^{-12}(\Sigma+1)$. Hence, if
$\hat s>\vartheta$,
\[
 s\ge\hat s-1.3\cdot10^{-15}\,\Sigma>0.9\cdot10^{-12}(\Sigma+1)-1.3\cdot10^{-15}\,\Sigma>0.
\]
The same conclusion holds with subnormal operands. Every exact input that
rounds to a finite binary64 number has absolute value less than
$2^{1025}$, so replacing a product of two exact inputs by the product of
their rounded values causes an additional error of at most
$2\eta(1+\uround)2^{1025}+\eta^2<2^{-48}$. With at most eight products and
nine arithmetic operations, an allowance of $10^{-13}$ covers all such
absolute errors. The bound for $|\hat s-s|$ becomes
$1.3\cdot10^{-15}\Sigma+10^{-13}$, while the lower bound
$\vartheta>0.9\cdot10^{-12}(\Sigma+1)$ still holds, because the
corresponding absolute errors in $\vartheta$ are multiplied by
$10^{-12}$. These bounds again give $s>0$.
\end{proof}

If $\hat s\le\vartheta$, the inequality $G_jx_I<g_j$ is checked exactly in
integers.

\begin{lemma}[distance filter]\label{lem:filter2}
Let $v,w\in\Q^n$, $n\le8$, have all coordinates less than $4$ in absolute
value, and let $\hat\delta$ be the value of
$|\hat v|^2+|\hat w|^2-2\ip{\hat v}{\hat w}$ computed in floating point
from the rounded coordinates. Then $\bigl|\hat\delta-|v-w|^2\bigr|<10^{-12}$.
\end{lemma}

\begin{proof}
All products $|v_tw_t|$, $v_t^2$, $w_t^2$ are less than $16$, so
$|v|^2$, $|w|^2$ and $|\ip vw|$ are less than $16n$. Rounding the
coordinates changes each of these three quantities by at most
$\gamma_2\cdot16n$, and the floating-point evaluation adds at most
$\gamma_n(1+\gamma_2)\cdot16n$; together this is at most
$\gamma_{n+2}\cdot16n<1.5\cdot10^{-13}$ for $n\le8$. Multiplication by $2$
is exact, and the two additions contribute at most
$\gamma_2\cdot64n(1+10^{-12})<1.2\cdot10^{-13}$. In total,
$\bigl|\hat\delta-|v-w|^2\bigr|<(1+1+2)\cdot1.5\cdot10^{-13}+1.2\cdot10^{-13}<10^{-12}$.
Since the coordinates are less than four in absolute value, gradual
underflow contributes less than $10^{-300}$, which is absorbed by the same
bound.
\end{proof}

The program computes $\hat\delta$ for all vertex pairs and recomputes in
rational arithmetic every pair with $\hat\delta>1-10^{-9}$. A pair with
$|v-w|^2\ge1$ has $\hat\delta>1-10^{-12}$ and is therefore always
recomputed, so the conclusion that all distances are less than one is
rigorous. The printed diameter is a floating-point approximation to
$\sqrt{\max\hat\delta}$, and the constants in Table~\ref{tab:contact}
allow an error of $2\cdot10^{-12}$ in the squared diameter, which covers
both the distance filter and the final rounding of the square root.
Whether a leaf is accepted depends only on the strict test, not on the
printed value.

\subsection{Checking a single cell}\label{app:cellcheck}

For cell $i$ of a box $Q$ the program proceeds as follows.
\begin{enumerate}[label=\arabic*.]
\item It assembles the system $Gx\le g$ from the $2n-7$ inequalities of the
frame $F_n$, the $k-1$ inequalities $\ip{a_j-a_i}x\le\beta_i-\beta_j$
defining $C_i$, and the tangent halfspaces $\Pi_{i,j}$.
\item For each tangent halfspace it checks the condition of
Lemma~\ref{lem:support} with the exact centre and squared radius of the
corresponding ball.
\item It finds an interior point of $K_i$ (the centre of a largest
inscribed ball, by linear programming) and uses Qhull to obtain the family
$\mathcal S$ of active sets of the putative vertices.
\item For each $I\in\mathcal S$ it solves $G_Ix=g_I$ in integers, checks
the solution and $\det G_I\not\equiv0\pmod{2^{61}-1}$, and checks the
strict inactivity of the remaining constraints (by
Lemma~\ref{lem:filter1}, and exactly where the filter is inconclusive).
\item It checks the hypothesis of Proposition~\ref{prop:closure}. If steps
3--5 fail (typically because nearly parallel halfspaces make the vertex
structure degenerate), it tries subsets of the halfspaces, one per ball and
per cluster of nearly equal normals; a subset of valid halfspaces defines a
larger polytope that still contains $\Omega_Q\cap C_i$, so success on a
subset suffices. As a last resort it enumerates the vertices with
\texttt{cddlib}.
\item It selects vertex pairs by the filter of Lemma~\ref{lem:filter2} and
checks $|v-w|^2<1$ exactly for the selected pairs.
\end{enumerate}
Floating point is used in steps 3, 4 and 6, but none of its results is
taken on trust. Qhull only proposes candidates, which are then verified
exactly, and a filter draws a conclusion by itself only when
Lemma~\ref{lem:filter1} or Lemma~\ref{lem:filter2} justifies it; otherwise
the constraint or the pair is passed on to exact recomputation.

\subsection{Runs and independent checks}\label{app:runs}

The complete verification runs, whose logs are contained in the
repository~\cite{Repo}, took about $75$ minutes for $n=4$ and about
$11.5$ hours for $n=5$ on four processes. (These are measurements of those
runs and only indicate the order of magnitude on other machines.) For both
dimensions the logs record that all leaves were checked, that there were
no failures, that the covered volume is one and that the tree is complete.

A second, independent program re-verified a random sample of $120$ leaves in each dimension. It
enumerates every vertex set by \texttt{cddlib} in GMP rational arithmetic,
detects unboundedness by the presence of rays, and uses neither
Proposition~\ref{prop:closure} nor subsets of halfspaces. These runs reported no
failures. They are an additional check of the vertex enumeration; the
proof for the whole parameter cube rests on the complete runs.

\subsection{Reproducing the computations}\label{app:repro}

The accompanying repository~\cite{Repo} contains the data for both
partition trees, the code that verifies them, the code that checks the
covering of Section~\ref{sec:six}, and the logs of the complete runs. The
two trees are stored in the files
\path{certificates/b4_k6/cert.json.gz} and
\path{certificates/b5_k12/cert.json.gz} (gzip-compressed JSON). The
preorder tree uses $0,1,2,3$ for a midpoint split in the corresponding
coordinate, $-(j+1)$ for a leaf with diagram $j$, and \texttt{null} for an
unresolved leaf; a complete tree has no unresolved leaf. Each diagram is a
pair $(A,b)$ of integer arrays with common denominator $10^9$. The
five-dimensional file records its dimension and number of parts
explicitly; for the four-dimensional file they are four and six. The
tangent halfspaces are not stored; they are regenerated and checked by the
verification program.

From the repository root, the commands
\begin{verbatim}
python3 verify.py inspect
python3 verify.py short
\end{verbatim}
use only the Python standard library. The first checks checksums, the tree
structure, the exact covered volume, the shapes of the diagrams and the
stored log summaries; it does not replace a geometric verification. The
second performs the checks of Appendix~\ref{app:rays} and the arithmetic
in the proofs of Theorems~\ref{thm:main678}--\ref{thm:ranges}. The complete geometric verifications are
run by
\begin{verbatim}
python3 verify.py full --dimension 4 --workers 4
python3 verify.py full --dimension 5 --workers 4
\end{verbatim}
and require NumPy, SciPy and pycddlib with GMP. Installation commands
and the optional independent \texttt{cddlib} checks are documented in
\texttt{README.md}.
The command \texttt{python3 verify.py all --workers 4} runs all verification
stages, including both complete geometric checks. Python assertions must
remain enabled, so neither \texttt{python -O} nor \texttt{PYTHONOPTIMIZE}
may be used.

The interval arithmetic derives from the
programs accompanying~\cite{ABPR}; this is recorded in the file
\texttt{THIRD\_PARTY\_NOTICES.md} of the repository.

\section{How the partitions were found}\label{app:search}

The procedures below produced the partition trees of
Section~\ref{sec:trees}. They play no role in the proof: everything they
output passes through the verification of Appendix~\ref{app:verify}. The
repository~\cite{Repo} contains the resulting certificates and the
programs needed for their verification.

\paragraph{Optimizing a partition.} For a fixed convex body $Y$ (a body
$H(u)$ or $\Omega_Q$) we minimize $\Phi=\max_i\diam(Y\cap C_i)$ over the
parameters $(a_i,\beta_i)$ of the power diagram. The body $Y$ is replaced by
a dense sample of boundary points, obtained by shooting rays from an
interior point. For each cell we select several pairs of sample points that
are farthest apart and start from each pair a local maximization of
$|x-y|^2$ subject to $x,y\in Y\cap C_i$ (SLSQP). By the envelope theorem,
the gradient of the local maximum with respect to $(a_i,\beta_i)$ is given
by the Lagrange multipliers of the active cell constraints. With these
gradients we make a linear programming step that minimizes the largest
linearized squared diameter within a trust region. The initial partition is
obtained by $k$-means clustering of the sample.

\paragraph{The tree.} We first optimized partitions at a few sample points
$u$ of the family, with $Y=H(u)$. The search then processes the tree from
the root. For a box it computes the uncertainty of the enclosure,
$\rho^+-\rho^-+\max_iE_i$, and splits the box at once if this number is
large. Otherwise it tries, in this order, the nearest sample partition,
transported to the box by the orthogonal map of
$\operatorname{span}(e_1,e_2)$ that best matches the two contact
triangles; the partition of the parent box; and a few optimization steps.
As soon as the numerical value of $\Phi$ drops below $0.99997$, an exact
verification of the leaf is attempted. A box that fails is split along the coordinate
whose halves have the smallest uncertainty.

\paragraph{Tangent halfspaces.} For cell $i$ we start with a few dozen
random halfspaces tangent to the four balls of $\Omega_Q$. We then compute
the vertices of $K_i$ in floating point and find the pairs at distance at
least one. For each vertex of such a pair that lies outside one of the
balls, we add the halfspace tangent to that ball at the projection of the
vertex. The normal is rounded to rationals with denominator $10^9$, and the
right-hand side is rounded up so that the condition of
Lemma~\ref{lem:support} holds exactly. This is repeated until all vertex
pairs are at distance less than one. The verifier runs the same procedure,
with a fixed seed of the random number generator for each cell, and then
checks every halfspace exactly.


\begin{thebibliography}{99}
\small

\bibitem{ABPR} A.~Arman, A.~Bondarenko, A.~Prymak and D.~Radchenko,
Borsuk number in four dimensions is at most seven, preprint (2026),
\href{https://arxiv.org/abs/2610.39078}{arXiv:2610.39078}; computational
supplement at \url{https://github.com/andriyprm/4d7balls}.

\bibitem{ABPRUniversal} A.~Arman, A.~Bondarenko, A.~Prymak and
D.~Radchenko, On asymptotic Lebesgue's universal covering problem,
\emph{Bull. Lond. Math. Soc.} \textbf{58} (2026), no.~10, e70514,
\href{https://doi.org/10.1112/blms.70514}{doi:10.1112/blms.70514}.

\bibitem{Aurenhammer} F.~Aurenhammer, Power diagrams: properties,
algorithms and applications, \emph{SIAM J. Comput.} \textbf{16} (1987),
78--96.

\bibitem{AvisFukuda} D.~Avis and K.~Fukuda, A pivoting algorithm for
convex hulls and vertex enumeration of arrangements and polyhedra,
\emph{Discrete Comput. Geom.} \textbf{8} (1992), 295--313.

\bibitem{Qhull} C.~B.~Barber, D.~P.~Dobkin and H.~Huhdanpaa, The
Quickhull algorithm for convex hulls, \emph{ACM Trans. Math. Software}
\textbf{22} (1996), 469--483.

\bibitem{Bareiss} E.~H.~Bareiss, Sylvester's identity and multistep
integer-preserving Gaussian elimination, \emph{Math. Comp.} \textbf{22}
(1968), 565--578.

\bibitem{BikeevRaigorodskii} A.~I.~Bikeev and A.~M.~Raigorodskii, On upper
bounds for the number of parts in the problem of partitioning sets into
parts of smaller diameter, \emph{Math. Notes} \textbf{119} (2026), no.~3,
406--415,
\href{https://doi.org/10.1134/S0001434626601899}{doi:10.1134/S0001434626601899}.

\bibitem{Bondarenko} A.~Bondarenko, On Borsuk's conjecture for two-distance
sets, \emph{Discrete Comput. Geom.} \textbf{51} (2014), 509--515.

\bibitem{Borsuk} K.~Borsuk, Drei S\"atze \"uber die $n$-dimensionale
euklidische Sph\"are, \emph{Fund. Math.} \textbf{20} (1933), 177--190.

\bibitem{BourgainLindenstrauss} J.~Bourgain and J.~Lindenstrauss, On
covering a set in $\R^N$ by balls of the same diameter, in \emph{Geometric
Aspects of Functional Analysis (1989--90)}, Lecture Notes in Math.~1469,
Springer, Berlin, 1991, 138--144.

\bibitem{Chalmers} L.~R.~Chalmers, The four-dimensional Borsuk number is
five or six, preprint, version~1, posted 4~November 2026,
\url{https://www.preprints.org/manuscript/202611.0339}.

\bibitem{DGK} L.~Danzer, B.~Gr\"unbaum and V.~Klee, Helly's theorem and its
relatives, in \emph{Convexity}, Proc. Sympos. Pure Math.~7, Amer. Math.
Soc., Providence, RI, 1963, 101--180.

\bibitem{Eggleston} H.~G.~Eggleston, Covering a three-dimensional set with
sets of smaller diameter, \emph{J. London Math. Soc.} \textbf{30} (1955),
11--24.

\bibitem{Filimonov} V.~P.~Filimonov, Covering planar sets, \emph{Sb. Math.}
\textbf{201} (2010), 1217--1248.

\bibitem{cdd} K.~Fukuda, \emph{cddlib}: an implementation of the double
description method, \url{https://github.com/cddlib/cddlib}.

\bibitem{FukudaProdon} K.~Fukuda and A.~Prodon, Double description method
revisited, in \emph{Combinatorics and Computer Science}, Lecture Notes in
Comput. Sci.~1120, Springer, Berlin, 1996, 91--111.

\bibitem{Grunbaum} B.~Gr\"unbaum, A simple proof of Borsuk's conjecture in
three dimensions, \emph{Proc. Cambridge Philos. Soc.} \textbf{53} (1957),
776--778.

\bibitem{HLP} G.~H.~Hardy, J.~E.~Littlewood and G.~P\'olya,
\emph{Inequalities}, 2nd ed., Cambridge University Press, 1952.

\bibitem{Heppes} A.~Heppes, T\'erbeli ponthalmazok feloszt\'asa kisebb
\'atm\'er\H oj\H u r\'eszhalmazok \"osszeg\'ere, \emph{Magyar Tud. Akad.
Mat. Fiz. Oszt. K\"ozl.} \textbf{7} (1957), 413--416.

\bibitem{Higham} N.~J.~Higham, \emph{Accuracy and Stability of Numerical
Algorithms}, 2nd ed., SIAM, Philadelphia, 2002.

\bibitem{HinrichsRichter} A.~Hinrichs and C.~Richter, New sets with large
Borsuk numbers, \emph{Discrete Math.} \textbf{270} (2003), 137--147.

\bibitem{JenrichBrouwer} T.~Jenrich and A.~E.~Brouwer, A $64$-dimensional
counterexample to Borsuk's conjecture, \emph{Electron. J. Combin.}
\textbf{21} (2014), no.~4, Paper~4.29.

\bibitem{Jung} H.~W.~E.~Jung, \"Uber die kleinste Kugel, die eine
r\"aumliche Figur einschlie\ss t, \emph{J. Reine Angew. Math.} \textbf{123}
(1901), 241--257.

\bibitem{KahnKalai} J.~Kahn and G.~Kalai, A counterexample to Borsuk's
conjecture, \emph{Bull. Amer. Math. Soc. (N.S.)} \textbf{29} (1993),
60--62.

\bibitem{Kalai} G.~Kalai, Some old and new problems in combinatorial
geometry I: around Borsuk's problem, in \emph{Surveys in Combinatorics
2015}, London Math. Soc. Lecture Note Ser.~424, Cambridge University
Press, 2015, 147--174.

\bibitem{Lassak} M.~Lassak, An estimate concerning Borsuk partition
problem, \emph{Bull. Acad. Polon. Sci. S\'er. Sci. Math.} \textbf{30}
(1982), 449--451.

\bibitem{Lenz} H.~Lenz, \"Uber die Bedeckung ebener Punktmengen durch solche
kleineren Durchmessers, \emph{Arch. Math.} \textbf{7} (1956), 34--40,
\href{https://doi.org/10.1007/BF01900521}{doi:10.1007/BF01900521}.

\bibitem{MooreKearfott} R.~E.~Moore, R.~B.~Kearfott and M.~J.~Cloud,
\emph{Introduction to Interval Analysis}, SIAM, Philadelphia, 2009.

\bibitem{Naszodi} M.~Nasz\'odi, Flavors of translative coverings, in
\emph{New Trends in Intuitive Geometry}, Bolyai Soc. Math. Stud.~27,
Springer, Berlin, 2018, 335--358.

\bibitem{Perkal} J.~Perkal, Sur la subdivision des ensembles en parties de
diam\`etre inf\'erieur, \emph{Colloq. Math.} \textbf{1} (1947), 45.

\bibitem{Raigorodskii} A.~M.~Raigorodskii, Borsuk's problem and the
chromatic numbers of some metric spaces, \emph{Russian Math. Surveys}
\textbf{56} (2001), 103--139.

\bibitem{RaigorodskiiLectures} A.~M.~Raigorodskii, Three lectures on the
Borsuk partition problem, in \emph{Surveys in Contemporary Mathematics},
London Math. Soc. Lecture Note Ser.~347, Cambridge University Press, 2007,
202--247.

\bibitem{Rockafellar} R.~T.~Rockafellar, \emph{Convex Analysis}, Princeton
University Press, Princeton, NJ, 1970.

\bibitem{Schramm} O.~Schramm, Illuminating sets of constant width,
\emph{Mathematika} \textbf{35} (1988), 180--189.

\bibitem{Repo} A.~Tolmachev, \emph{borsuk-bounds}: certificates and verification
programs for this paper, GitHub repository (2026), \repo.

\bibitem{TPV} A.~D.~Tolmachev, D.~S.~Protasov and V.~A.~Voronov, Coverings
of planar and three-dimensional sets with subsets of smaller diameter,
\emph{Discrete Appl. Math.} \textbf{320} (2022), 270--281.

\bibitem{TV} A.~Tolmachev and V.~Voronov, Reducing the upper bound for the
Borsuk number in $\R^4$ to $8$, preprint (2026),
\href{https://arxiv.org/abs/2605.19068}{arXiv:2605.19068}.

\bibitem{Zong} C.~Zong, Borsuk's partition conjecture, \emph{Jpn. J. Math.}
\textbf{16} (2021), 185--201.

\end{thebibliography}
\end{document}